\documentclass{jfm}

\usepackage{graphicx}
\usepackage{natbib}
\usepackage{amsmath}
\usepackage[colorlinks,pdftex]{hyperref}  

\ifCUPmtlplainloaded \else
 \checkfont{eurm10}
 \iffontfound
   \IfFileExists{upmath.sty}
     {\typeout{^^JFound AMS Euler Roman fonts on the system,
                  using the 'upmath' package.^^J}%
      \usepackage{upmath}}
     {\typeout{^^JFound AMS Euler Roman fonts on the system, but you
                  dont seem to have the}%
      \typeout{'upmath' package installed. JFM.cls can take advantage
                of these fonts,^^Jif you use 'upmath' package.^^J}%
     }
 \else
 \fi
\fi

\ifCUPmtlplainloaded \else
 \checkfont{msam10}
 \iffontfound
   \IfFileExists{amssymb.sty}
     {\typeout{^^JFound AMS Symbol fonts on the system, using the
               'amssymb' package.^^J}%
      \usepackage{amssymb}%

     }{}
 \fi
\fi

\ifCUPmtlplainloaded \else
 \IfFileExists{amsbsy.sty}
   {\typeout{^^JFound the 'amsbsy' package on the system, using it.^^J}%
    \usepackage{amsbsy}}
   {\providecommand\boldsymbol[1]{\mbox{\boldmath $##1$}}}
\fi

\providecommand\bnabla{\boldsymbol{\nabla}}
\providecommand\bcdot{\boldsymbol{\cdot}}

\newsavebox{\astrutbox}
\sbox{\astrutbox}{\rule[-5pt]{0pt}{20pt}}

\newcommand{\reals}{{\mathbb R}}
\newcommand{\mbs}[1]{\boldsymbol{#1}}

\title[The role of membrane viscosity in the dynamics of fluid membranes]{The role of membrane viscosity in the dynamics of fluid membranes}
\author[M.~Arroyo, A.~DeSimone and L.~Heltai]%
{M\ls A\ls R\ls I\ls N\ls O \ns A\ls R\ls R\ls O\ls Y\ls O,$^1$%
\ns A\ls N\ls T\ls O\ls N\ls I\ls O \ns D\ls E\ls S\ls I\ls M\ls O\ls N\ls E$^2$ \ns \and L\ls U\ls C\ls A \ns  H\ls E\ls L\ls T\ls A\ls I$^2$}

\affiliation{$^1$Dept. Applied Mathematics 3, LaC\`aN, Universitat Polit\`ecnica de Catalunya (UPC), Barcelona 08034, Spain\\[\affilskip]
$^2$SISSA-International School for Advanced Studies, Via Bonomea 265,
34136 Trieste, Italy}

\pubyear{2008}
\volume{???}
\pagerange{???--???}
\date{\today}

\begin{document}

\maketitle

\begin{abstract}
Fluid membranes made out of lipid bilayers are the fundamental separation structure in eukaryotic cells. Many physiological processes rely on dramatic shape and topological changes (e.g. fusion, fission) of fluid membrane systems. Fluidity is key to the versatility and constant reorganization of lipid bilayers. Here, we study the role of the membrane intrinsic viscosity, arising from the friction of the lipid molecules as they rearrange to accommodate shape changes, in the dynamics of morphological changes of fluid vesicles. In particular, we analyze the competition between the membrane viscosity and the viscosity of the bulk fluid surrounding the vesicle as the dominant dissipative mechanism. We consider the relaxation dynamics of fluid vesicles put in an out-of-equilibrium state, but conclusions can be drawn regarding the kinetics or power consumption in regulated shape changes in the cell. On the basis of numerical calculations, we find that the dynamics arising from the membrane viscosity are qualitatively different from the dynamics arising from the bulk viscosity. When these two dissipation mechanisms are put in competition, we find that for small vesicles the membrane dissipation dominates, with a relaxation time that scales as the size of the vesicle to the power 2. For large vesicles, the bulk dissipation dominates, and the exponent in the relaxation time vs. size relation is 3. 
\end{abstract}

\section{Introduction}

Amphiphilic membranes are self-assembled structures made out of lipids or other amphiphilic molecules such as diblock co-polymers \citep{Discher:1999sp}. Above a transition temperature, these membranes are fluid within the membrane surface, while they retain the transversal order, which confers them with bending rigidity. The membrane fluidity is essential for many biochemical processes involving membrane proteins \citep{Saffman08011975}. The fluidity of amphiphilic membranes has been studied experimentally, with measurements of its two-dimensional viscosity \citep{DanovDimova00}, as well as computationally with molecular dynamics (MD) simulations of viscometric tests \citep{Otter-Biophis-J}. Atomistic molecular dynamics (MD) simulations have been very useful in this and other respects, but remain limited to small membrane patches due to the prohibitive number of atoms involved in closed vesicles and the slow equilibration times of membrane systems \citep{ISI:000230886900046}. Coarse-grained MD allows us to reach larger systems, and there has been notable successful studies in recent years involving out-of-equilibrium phenomena at the scale of small vesicles, e.g. \cite{ISI:000246693100044}. Yet, as in atomistic MD, the computational cost scales as the size of the system to the power 6, and therefore there exists a hard upper bound on the sizes of systems that can be reached with current computers. On the other end of the spectrum of modeling approaches for fluid membranes, continuum mechanics has proven very effective in describing the mechanics of membrane systems and reproduce experimentally observed equilibrium shapes \citep{Lipowsky95}. There has been much less progress in using continuum mechanics models to describe the dynamics of fluid membrane systems, and in particular to account for the membrane viscous flow, i.e. the rearrangement of lipids on the membrane surface driven by external actions or required to accommodate shape changes. This scarcity of results contrasts with the importance of out-of-equilibrium and dynamical events both in the cell and in synthetic systems \citep{Baumgart03,Bacia:2005lr}. Here, we address the relaxation dynamics of vesicles (closed membranes made out of fluid amphiphilic membranes) with continuum mechanics models and simulations, with a particular emphasis on the role of the membrane dissipation due to the friction between the amphiphiles as they shear to accommodate shape changes, and the competition in setting the relaxation dynamics between this membrane flow and the flow induced in the ambient fluid by shape changes. Our simulations can treat without difficulties large systems and very slow processes, in contrast with MD methods, at the expense of the molecular details, which are crucial for some phenomena.

The continuum mechanics formulation of the coupling between shape dynamics and membrane hydrodynamics was first established by \cite{Scriven60} (see also \cite{Aris89}). The Navier-Stokes equations on a two-dimensional, time-evolving, curved manifold were given in an intrinsic manner in the parameter space of the surface, and required heavy differential geometry artillery. We note that the inertial terms in the surface Navier-Stokes equations are easily dealt with, and the main complications arise from the viscous terms; here we are interested in the low Reynolds number regime, and therefore ignore the inertial terms. An alternative cartesian formulation in terms of time-dependent projection operators was proposed later \citep{Secomb82,BarthesJFM}. None of these formulations have been exercised beyond surface flows on fixed, simple shapes, or beyond infinitesimal shape perturbations around simple shapes such as the sphere. For finite shape changes, only extremely restricted families of idealized shapes have been considered, e.g. \cite{Fischer94}. Maybe due to the inherent complexity of the equations governing the coupled shape dynamics and membrane flow, the effect of membrane viscosity has often been neglected in continuum studies, and the idea that this source of dissipation can be safely ignored in most situations of interest, compared for instance to bulk dissipation, has prevailed. Here, we challenge this viewpoint, which of course can be well justified in important situations such as the behavior of vesicles in shear flow \citep{pozrikidis,misbahPRL07,Misbah_PRE_05}. But even in this case, recent coarse-grained MD simulations suggest that the membrane viscosity can play a crucial role \citep{ISI:000230886900072,ISI:000232392900009}.

We present, to the best of our knowledge, the first vesicle calculations involving finite shape changes and incorporating the effect of the membrane flow (as well as the effect of the bulk flow). We perform relaxation dynamics simulations of axisymmetric vesicles driven by curvature elasticity, dragged by membrane and bulk viscosity, and subject to the usual incompressibility and inextensibility constraints of the bulk and the membrane fluids respectively. The new geometric and variational formulation of the Stokes equations on a curved, time-evolving surface coupled to the ambient fluid flow and the curvature elasticity proposed in \cite{arroyo-desimone-2009} is the basis of our numerical strategy.

The mathematical formulation of the relaxation dynamics of axisymmetric vesicles is given in Section \ref{formulation}. The general equations are particularized, and strikingly simple expressions are developed for the membrane dissipation functional and the local inextensibility constraint, which make the theory accessible to analytical calculations and simple numerical implementations. The numerical discretization of the theory, with Galerkin finite elements and boundary elements based on B-Splines is described in Section \ref{numerical}. The numerical results investigating the qualitative behavior of the dynamics engendered by different dissipation mechanisms, and exploring the competition between membrane and bulk viscosity are reported in Section \ref{results}. The conlcusions we can draw from our analysis are collected in Section \ref{conclusions}.

\section{Formulation of the relaxation dynamics}
\label{formulation}

This section describes the continuum model for an inextensible fluid membrane with curvature elasticity immersed in a viscous, incompressible Newtonian fluid. We ignore inertial forces, and model the amphiphilic membrane as a Newtonian two-dimensional fluid, in agreement with coarse-grained molecular dynamics simulations \citep{Otter-Biophis-J} and experimental observations \citep{DanovDimova00,Dimova:2006ys,Cicuta:2007lr}. Besides the effect of the membrane viscosity, which is original to the present work, the ingredients in the formulation are classical and well-known, although we introduce some particularly simple expressions for axisymmetric vesicles. We firstly describe the kinematics, which we choose to describe in a Lagrangian manner both in the normal and tangential directions. We then describe the dissipative viscous mechanisms in terms of the Rayleigh dissipation potentials, the energetic mechanism operative in this model (curvature elasticity), and the constraints. We derive the force balance equations from a variational principle, from which the relaxation dynamics follow.

\subsection{Kinematics}

We describe parametrically axisymmetric vesicles in terms of the generating curve, i.e. the vesicle surface $\Gamma_t$ at a given instant $t$ is given by 
\[
\mbs{x}(u,\theta; t) = \{r(u;t)\cos\theta, r(u;t)\sin\theta,z(u;t) \}, ~~~~u\in [0,1], ~\theta\in [0,2\pi],
\]
where
\[
\mbs{c}(u; t) = \{r(u; t), z(u; t) \}, ~~~~u\in [0,1],
\]
is the parametric description of the generating curve $\mathcal{C}_t$ at the instant $t$. We consider closed surfaces with continuous tangents, hence we require 
\begin{equation}
r(0) = r(1) = 0, ~~~~~ z'(0) = z'(1)= 0,
\label{BCs}
\end{equation}
where $(\cdot)'$ denotes partial differentiation with respect to $u$.

For simplicity, from this point on we omit the dependence of all quantities on time. We shall formulate the mechanics of the fluid membrane in terms of the generating curve. Its speed is given by $a(u) = \sqrt{[r'(u)]^2+[z'(u)]^2}$. Integrals on the surface can be brought to the interval $[0,1]$ with the relation ${\rm d}S = (2\pi a r){\rm d}u$. The tangent unit vector to the generating curve pointing in the $u$ direction and a unit normal are given by
\[
\mbs{t}= \frac{1}{a}\{r',z'\}, ~~~~~~~\mbs{n}= \frac{1}{a}\{-z',r'\}.
\]
At this point, there are several possibilities in describing the kinematics of the membrane. Since amphiphilic membranes are two-dimensional fluid bodies moving in a higher-dimensional space, the description of the motion is necessarily Lagrangian, at least with regards to the shape of the membrane, i.e. $\dot{\mbs{c}} \bcdot \mbs{n} = v_n$ on $\Gamma$, where the dot denotes partial differentiation with respect to time and $v_n$ is the normal velocity of the membrane. In the tangential direction, however, there is gauge freedom due to the particle relabeling symmetry of the fluid flow equations, here the re-parameterization invariance with respect to $u$. One special choice is the Eulerian gauge, for which $\dot{\mbs{c}} \bcdot \mbs{t} = 0$ and $v_t$, the tangential velocity of the material particles on the membrane, enters the formulation as an independent variable. Here we adopt another natural choice, the Lagrangian gauge, in which the parameter $u$ is viewed as a label for material particles, hence $\dot{\mbs{c}} \bcdot \mbs{t} = v_t$. In this case, we have:
\[
\dot{\mbs{c}} = \{ \dot{r}, \dot{z}\} = \mbs{V} = \mbs{v} + v_n \mbs{n}= v_t \mbs{t} + v_n \mbs{n},
\]
and
\begin{equation}
v_t = \frac{1}{a}(r'\dot{r} + z'\dot{z}), ~~~~ v_n = \frac{1}{a}(-z'\dot{r} + r'\dot{z}).
\label{v_rz}
\end{equation}
This decomposition is important since only the normal velocities change the shape of the membrane, while the tangential velocities represent the flow of the amphiphiles on the membrane surface. The above relations can be inverted as
\[
\dot{r} = \frac{1}{a}(r'v_t - z'v_n), ~~~~ \dot{z} = \frac{1}{a}(r'v_n + z'v_t).
\]
A key object for the kinematics and constitutive relations of the two-dimensional fluid on a curved, time-evolving surface is the rate-of-deformation tensor $\mbs{d}$. This tensor can be geometrically thought of as the tangent projection of the rate of change of the metric tensor of the surface as it is advected by the membrane velocity  \citep{Scriven60,MarsdenHughes,arroyo-desimone-2009}, which leads to the expression 
\[
\mbs{d} = \frac{1}{2}\left(\bnabla_s \mbs{v} + \bnabla_s \mbs{v}^T\right) - v_n \mbs{k},
\]
where $\bnabla_s$ denotes the surface covariant derivative and $\mbs{k} = -\bnabla_s\mbs{n}$ the second fundamental form. The formulas for axisymmetric surfaces can be found in \cite{arroyo-desimone-2009}.

Introducing $b(u) = -r''(u)z'(u) + r'(u) z''(u)$, we can write the mean and gaussian curvatures as
\begin{equation}
H = \frac{1}{a}\left(\frac{b}{a^2} + \frac{z'}{r} \right), ~~~~ K = \frac{b z'}{a^4 r}.
\label{curvatures}
\end{equation}
Note that here $H$ is defined as the trace of the second fundamental form.

\subsection{Governing equations}

We present first the generic ingredients entering the governing equations for the shape evolution, under the assumption of low Reynolds number, hence neglecting inertia. The Rayleigh dissipation potentials can be expressed by a functional $W[\dot{r},\dot{z}]$, which, although not explicitly written, depends nonlinearly on the shape of the surface. Under the assumption that both the ambient bulk fluid and the membrane two-dimensional fluid are Newtonian, this functional is quadratic in its arguments. 

The curvature energy and in general other energetic mechanisms such as line tension, depend exclusively on the shape of the membrane, hence can be expressed in terms of a nonlinear functional $\Pi[r,z]$. The energy release rate functional is the negative of the variation of the energy functional in the direction of $\{\dot{r},\dot{z} \}$, i.e. $G[\dot{r},\dot{z}] = - \delta\Pi[r,z; \dot{r},\dot{z}]$, where again its explicit (nonlinear) dependence on the configuration of the surface has been omitted. By its definition, it is clear that the energy release rate functional is linear in its arguments.

The membrane dynamics are often constrained by global nonlinear equalities such as enclosed volume or surface area constraints. We generically express these constraints with a vector-valued functional as $\mbs{B}[r,z]=\mbs{0}$, or in rate form as $\mbs{C}[\dot{r},\dot{z}] = \delta \mbs{B}[r,z; \dot{r},\dot{z}] = \mbs{0}$. Again, this functional is linear in its arguments and nonlinear in the configuration of the membrane. We shall also consider local constraints such as the local inextensibility of the membrane, expressed at each point of the membrane as $c(\dot{r},\dot{z})= 0$.

At each instant, the dynamics of the membrane equilibrate the dissipative and the energetic forces subject to the constraints. In other words, the membrane evolves in such a way that the rate of viscous dissipation resisting the flow is exactly equal to the energy release rate during the flow. Mathematically, the evolution equations follow from minimizing $W[\dot{r},\dot{z}]-G[\dot{r},\dot{z}]$ subject to the constraints. Forming the Lagrangian
\[
\mathcal{L}[\dot{r},\dot{z},\lambda, \mbs{\Lambda}] = W[\dot{r},\dot{z}]-G[\dot{r},\dot{z}]  - \int_{\Gamma} \lambda c(\dot{r},\dot{z}) {\rm d}S - \mbs{\Lambda}\bcdot \mbs{C}[\dot{r},\dot{z}],
\]
the velocities and the Lagrange multipliers at each configuration $\{r,z\}$ are obtained as the stationary points with respect to all admissible variations
\begin{equation}
\delta_{\dot{r}}\mathcal{L} = \delta_{\dot{z}}\mathcal{L} = \delta_{\lambda}\mathcal{L} = \delta_{\mbs{\Lambda}}\mathcal{L} = 0,
\label{pvp}
\end{equation}
which is a form of the principle of virtual power. Note that, unlike the other arguments of the Lagrangian functional, the Lagrange multipliers for the global constraints $\mbs{\Lambda}$ are not functions of $u$. If $c(\dot{r},\dot{z})= 0$ expresses the local inextensibility of the membrane, then the Lagrange multiplier $\lambda$ is the surface tension. We provide below the specific form of these functionals.

\subsection{Constraints}

Fluid membranes are semipermeable, and assuming osmotic equilibrium between the enclosed and the outer media, it is often reasonable to assume that the enclosed volume is constant \citep{Handbook-Biological-Physics}. This constraint is expressed as
\begin{equation}
0 = C^{\rm vol}[\dot{r},\dot{z}] = \dot{V} = -\int_\Gamma v_n {\rm d}S =  -\int_0^1 \frac{1}{a}(-z'\dot{r} + r'\dot{z}) (2\pi a r){\rm d}u.
\label{D_vol} 
\end{equation}
Generally, the in-plane stresses on the membrane are low compared to the elastic moduli, and the membrane can be assumed to be locally inextensible \citep{Handbook-Biological-Physics}. This condition on the surface can be expressed as \citep{arroyo-desimone-2009}
\[
0 = {\rm trace}~\mbs{d} = \bnabla_s\bcdot \mbs{v}-H v_n,
\]
which for axisymmetric surfaces reduces to
\begin{equation}
0 = c^{\rm inext}(\dot{r},\dot{z}) = \frac{1}{ar} (r v_t)' - H v_n.
\label{inext_axi} 
\end{equation}
With the Lagrangian gauge mentioned earlier, this constraint takes the particularly simple expression (see appendix \ref{inext} for a derivation)
\begin{equation}
0 = c^{\rm inext}(\dot{r},\dot{z}) = \frac{\partial}{\partial t} \ln ar = \frac{\dot{a}}{a} + \frac{\dot{r}}{r} = \frac{1}{a^2} (r' \dot{r}' + z' \dot{z}') + \frac{\dot{r}}{r}.
\label{inext_axi_lagr} 
\end{equation}
If the details of the fluid flow on the membrane are disregarded, this local condition is often replaced by a total surface area constraint \citep{Du-Flow}, which for a closed surface takes the form
\begin{equation}
0 = C^{\rm area}[\dot{r},\dot{z}] = \dot{S} = \int_\Gamma {\rm trace}~\mbs{d} ~{\rm d}S = -\int_\Gamma H v_n {\rm d}S =  -\int_0^1 \frac{H}{a}(-z'\dot{r} + r'\dot{z}) (2\pi a r){\rm d}u,
\label{D_area} 
\end{equation}
where the divergence theorem has been used.

\subsection{Membrane dissipation}

As derived in \cite{arroyo-desimone-2009}, the Rayleigh dissipation potential for a Newtonian closed fluid membrane can be written as
\begin{eqnarray*}
W^{\rm mem}[\mbs{v},v_n] & = & \int_\Gamma \mu ~\mbs{d}\mbs{:}\mbs{d}~{\rm d}S \\ & = &
\int_\Gamma \mu \left[ \frac{1}{2} \vert {\bf d}\mbs{v}^\flat\vert^2 + (\bnabla_s\bcdot \mbs{v})^2 - K \vert \mbs{v} \vert^2 + (H^2-2K) v_n^2 - 2(\bnabla_s \mbs{v}\mbs{:}\mbs{k}) v_n\right] {\rm d}S,
\end{eqnarray*}
where ${\bf d}$ denotes the exterior derivative, ${\bf d}\mbs{v}^\flat$ is the generalization of the curl of the tangent velocity field on the surface, and $\mu$ denotes the membrane viscosity, with units of force $\times$ time $\times$ length$^{-1}$. Molecular dynamics simulations \citep{Otter-Biophis-J} as well as experiments \citep{DanovDimova00,Dimova:2006ys} support modeling fluid membranes as Newtonian two-dimensional fluids. For axisymmetric surfaces, this expression reduces to
\begin{eqnarray*}
W^{\rm mem}[\dot{r},\dot{z}] & = & \int_{\Gamma} \mu \left[ \left(\frac{1}{a} v_t' \right)^2 +  \left(\frac{r'}{a r} v_t \right)^2 - \frac{2}{a} v_n \left(\frac{b}{a^3} v_t' + \frac{z'r'}{a r^2} v_t  \right)  +  (H^2-2K) v_n^2  \right] {\rm d}S, \\
& = & \int_{\Gamma} \mu \left[\begin{array}{ccc}v_t' &  v_t & v_n\end{array}\right]
\left[\begin{array}{ccc}\frac{1}{a^2} & 0 & -\frac{b}{a^4} \\0 & \left(\frac{r'}{ar}\right)^2 & -\frac{z'r'}{(ar)^2} \\ -\frac{b}{a^4} & -\frac{z'r'}{(ar)^2} & H^2-2K \end{array}\right]\left[\begin{array}{c}v_t' \\  v_t \\ v_n\end{array}\right]
{\rm d}S.
\end{eqnarray*}
The matrix in the above expression can be checked to be positive semidefinite with rank 2, which reflects the zero-dissipation mode resulting from translating rigidly the vesicle along the symmetry axis in the $z$ direction, as a consequence of the internal nature of the membrane dissipation.

We provide here a remarkably simple alternative expression of this dissipation potential under the hypothesis of axisymmetry with the Lagrangian gauge. We first use the the local inextensibility constraint to express $v_t'$ in terms of $v_t$ and $v_n$. The dissipation potential then takes the form
\[
W^{\rm mem}[\dot{r},\dot{z}]  = 
\int_{\Gamma} \mu \frac{2}{(a r)^2} \left[\begin{array}{cc}  v_t & v_n\end{array}\right]
\left[\begin{array}{cc} (r')^2 & -r' z'  \\ -r' z' & (z')^2 \end{array}\right]\left[\begin{array}{c} v_t \\ v_n\end{array}\right]
{\rm d}S,
\]
which, using equation (\ref{v_rz}), reduces to 
\begin{equation}
W^{\rm mem}[\dot{r},\dot{z}]  = \int_{\Gamma} 2\mu \left( \frac{\dot{r}}{r} \right)^2 {\rm d}S = \int_0^1 2\mu \left( \frac{\dot{r}}{r} \right)^2   (2\pi a r){\rm d}u.
\label{Wmem}
\end{equation}

\subsection{Bulk dissipation}

We assume that the velocities of the membrane and of the surrounding fluid coincide, i.e. no slip as suggested by coarse-grained MD simulation \citep{Otter-Biophis-J} and conventionally assumed for lipid-water interactions \citep{Stone98}. This hypothesis may not be adequate in extreme situations, for instance if there is significant flow across the membrane. We consider an infinite incompressible fluid at rest at infinity and axisymmetry.

%

In the absence of body forces, the Rayleigh potential accounting for the dissipation in the surrounding fluid can be transformed into a surface integral by virtue of the divergence theorem as
\[
W^{\rm bulk} = \frac{1}{2} \int_{\reals^3} \mbs{\sigma} \mbs{:} \mathsfbi{D} ~{\rm d}V = \frac{1}{2} \int_\Gamma \Delta\mbs{g}\bcdot \mbs{V} {\rm d}S,
\]
where $\mbs{\sigma}$ is the stress tensor and $\mathsfbi{D}$ the rate-of-deformation tensor in the bulk fluid, $\mbs{V}$ denotes the velocity field on the membrane, and $\Delta\mbs{g}$ denotes the jump of tractions across the membrane due to the elastic forces, constraint reactions, and membrane viscous forces. By inverting the boundary integral representation of the velocities on the membrane in terms of a single layer potential, which assumes that the fluid viscosity is the same inside and outside of the vesicle, we obtain
\begin{equation}
\mbs{V}(u) = - \frac{1}{\mu^{\rm bulk}}  \int_\Gamma
\mathsfb{G}(\mbs{c}(u), \mbs{c}(v)) \bcdot  \Delta\mbs{g}(v)  ~{\rm
 d}S(v) := \frac{1}{\mu^{\rm bulk}} (\mathcal{A} \Delta\mbs{g})(u),
\label{bir}
\end{equation}
where $\mathsfb{G}$ is the axisymmetric Green's function and $\mu^{\rm bulk}$ is the bulk fluid viscosity. This representation allows us to implicitly express $W^{\rm bulk}$ in terms of $\dot{r}$ and $\dot{z}$.

The axisymmetric kernel $\mathsfb{G}$ represents the single layer
integral of a \emph{ring} of singularities passing through the point
$\mbs{c}(u)$:
\begin{equation}
 \label{axisymmetric-stokes-kernel}
 \mathsfb{G} (\mbs{c}(u), \mbs{c}(v)) :=  r(v) \int_0^{2\pi}
 \left[
   \begin{array}{cc}
     \mathsfb{S}_{zz} & (\mathsfb{S}_{yz} \cos\theta + \mathsfb{S}_{xz} \sin\theta) \\
     \mathsfb{S}_{zy} & (\mathsfb{S}_{yy} \cos\theta + \mathsfb{S}_{xy} \sin\theta)
   \end{array}
 \right]  ~{\rm d} \theta,
\end{equation}
where $\mathsfb{S}$ is the three-dimensional \emph{Stokeslet} 
\begin{equation*}
\mathsfb{S}(\mbs{r}) := \frac{1}{8\pi}\left(\frac{\mbs{I}}{|\mbs{r}|} +\frac{\mbs{r}\otimes\mbs{r}}{|\mbs{r}|^3}\right),
\end{equation*}
evaluated at $\mbs{r} = \mbs{R}(\theta) \mbs{c}(v) - \mbs{c}(u)$, with
$\mbs{R}(\theta)$ a rotation of angle $\theta$ around the
axis of symmetry.

A detailed derivation of equation~\eqref{axisymmetric-stokes-kernel}
can be found in~\cite{pozrikidis-book-1992}, together with an explicit
expression of $\mathsfb{G}$ in terms of complete elliptic integrals of
the first and second kind. Using the implicit definition of the
Dirichlet to Neumann map $\mathcal A^{-1}$ given in equation~\eqref{bir},
we can then express the bulk dissipation as
\begin{equation}
 \label{W-bulk}
 W^{\rm bulk}[\dot{r},\dot{z}]  = \int_{\Gamma} \frac{1}{2} \mu^{\rm bulk}
 \mbs{V}\cdot \mathcal A^{-1} \mbs{V}{\rm d}S = \int_0^1 \frac{1}{2}\mu^{\rm bulk}
 \{\dot{r}, \dot{z}\}\cdot \mathcal A^{-1} \{\dot{r}, \dot{z}\} (2\pi a r){\rm d}u.
\end{equation}
This expression is valid under the compatibility condition $\int_\Gamma \mbs{V}\bcdot
 \mbs{n} ~{\rm d}S = 0$, a result of the fluid incompressibility, enforced by Eq.~(\ref{D_vol})

\subsection{Curvature elasticity}

The curvature elasticity of the fluid membrane is modeled with the Helfrich-Canham functional \cite[see][for a discussion on curvature elasticity models]{Handbook-Biological-Physics},
\begin{equation}
\Pi = \int_{\Gamma} \frac{\kappa}{2} (H-C_0)^2~dS + \int_{\Gamma} \kappa_G K~{\rm d}S,
\end{equation}
where $C_0$ denotes the spontaneous curvature and $\kappa$ and $\kappa_G$ are elastic parameters. For closed homogeneous surfaces, the second term is a topological invariant, which we will disregard here. For a more general treatment including membranes with boundary, see \cite{arroyo-desimone-2009}. The energy release rate takes the form
\[
G[v_n] = -\int_\Gamma \kappa (H-C_0)\left[ \Delta_s v_n + \frac{1}{2}\left( H^2 - 4 K + H C_0  \right) v_n \right] {\rm d}S.
\]
For the numerical implementation, to avoid third order derivatives of $r$ and $z$, it proves more convenient to take variations directly from the following form of the elastic energy
\begin{equation}
\Pi[r,z] = \int_0^1 \frac{\kappa}{2} (H-C_0)^2 (2\pi a r){\rm d} u,
\label{curv_e}
\end{equation}
together with equation (\ref{curvatures}), which yields an expression for $G[\dot{r},\dot{z}]$ given in Appendix \ref{G}. Both approaches can be shown to be equivalent by integration by parts with the conditions in equation (\ref{BCs}) to annihilate the boundary terms.

\subsection{Different models}

In the paper, we consider several models to assess the specific effects of the membrane viscous flow and its relevance. Initially, for the sake of qualitative comparison, we consider a model with only membrane dissipation, and local membrane inextensibility
\[
\mathcal{L}^A[\dot{r}, \dot{z},\lambda,\Lambda^{\rm vol}] = W^{\rm mem}[\dot{r}, \dot{z}] - G[\dot{r}, \dot{z}] - \int_\Gamma \lambda c^{\rm inext}(\dot{r}, \dot{z}) {\rm d}S - \Lambda^{\rm vol} C^{\rm vol}[\dot{r}, \dot{z}].
\]
We consider also a model with only bulk dissipation, together with a global area constraint as in \cite{Du-Flow}
\[
\mathcal{L}^B[\dot{r}, \dot{z},\Lambda^{\rm area},\Lambda^{\rm vol}] = W^{\rm bulk}[\dot{r}, \dot{z}] - G[\dot{r}, \dot{z}] - \Lambda^{\rm area} C^{\rm area}[\dot{r}, \dot{z}] - \Lambda^{\rm vol} C^{\rm vol}[\dot{r}, \dot{z}].
\]
A similar model with local inextensibility can be considered, as in \cite{Misbah_PRE_05}. We have checked that the results for both models are very close, hence consider only $\mathcal{L}^B$ for definiteness. 

We also consider a model whose dissipative mechanism is of purely mathematical intent
\[
\mathcal{L}^C[\dot{r}, \dot{z},\Lambda^{\rm area},\Lambda^{\rm vol}] = W^{L_2}[\dot{r}, \dot{z}] - G[\dot{r}, \dot{z}] - \Lambda^{\rm area} C^{\rm area}[\dot{r}, \dot{z}] - \Lambda^{\rm vol} C^{\rm vol}[\dot{r}, \dot{z}],
\]
where
\[
W^{L_2}[\dot{r}, \dot{z}] = \int_\Gamma \frac{\hat{\mu}}{2}v_n^2{\rm d}S,
\]
and $\hat{\mu}$ is a mathematical viscosity coefficient. The resulting dynamics, for the closed surfaces considered here, are a constrained version of the $L_2$ gradient flow of the Willmore energy. Such a gradient flow finds applications in geometrical analysis \citep{MR1827100}, and is also considered by some authors as a simple model for the dynamics of fluid membranes \citep{Du-etal-jcp}. In all the models presented this far, vesicles of different sizes evolve in the exact same way, upon re-scaling of the time variable. Finally, the following model accounts for both the membrane and the bulk dissipation
\[
\mathcal{L}^{\rm full}[\dot{r}, \dot{z},\lambda,\Lambda^{\rm vol}] = W^{\rm mem}[\dot{r}, \dot{z}] + W^{\rm bulk}[\dot{r}, \dot{z}]  - G[\dot{r}, \dot{z}] - \int_\Gamma \lambda d^{\rm inext}(\dot{r}, \dot{z}) {\rm d}S - \Lambda^{\rm vol} D^{\rm vol}[\dot{r}, \dot{z}].
\]
In this model, the two dissipative mechanisms compete, and vesicles of different sizes exhibit a different relaxation behavior, as reported in \cite{arroyo-desimone-2009} and illustrated later in the present paper.

These models could be supplemented by external actions, although here we only consider the relaxation dynamics from an initial out-of-equilibrium condition. The out-of-equilibrium configurations are obtained as equilibria for a given value of the spontaneous curvature $C_0$, which in the dynamics simulation is set to zero, or by applying an external force on the vesicle, which is suddenly released so that the systems returns to an equilibrium state. It should be mentioned that, while for models $\mathcal{L}^B$, $\mathcal{L}^C$, and $\mathcal{L}^{\rm full}$ the boundary conditions in equation (\ref{BCs}) are sufficient, the model $\mathcal{L}^A$ requires an additional condition to fix the motion along $z$, which otherwise remains undetermined due to the invariance of $W^{\rm mem}$ apparent from equation (\ref{Wmem}).

\section{Numerical approximation}
\label{numerical}

\subsection{Spatial semi-discretization}

The spatial discretization follows from a standard Galerkin approach. The generating curve of the axisymmetric surface is represented numerically as a B-Spline curve
\[
\mbs{c}(u;t)= \{r(u; t), z(u; t) \} \approx \sum_{I=1}^N B_I(u)\underbrace{\{r_I(t), z_I(t) \}}_{\mbs{P}_I(t)},
\]
where $B_I(u)$ are the B-Spline basis functions \citep{nurbs} defined on the interval $[0, 1]$, and $\{r_I(t), z_I(t) \}$ is the position of the $I-$th control point of the B-Spline curve at instant $t$. Again, we drop the dependence on time. The velocity of the membrane can be computed as
\begin{equation}
\mbs{V}(u) \approx \sum_{I=1}^N B_I(u) \dot{\mbs{P}}_I.
\label{vh}
\end{equation}
Since the formulation involves up to second derivatives of the generating curve, it is convenient that the numerical representation is sufficiently smooth, hence avoiding cumbersome mixed approaches. Here, we have considered cubic B-Splines, which have up to second-order continuous derivatives. With the natural numbering of the basis functions, the symmetry conditions in equation (\ref{BCs}) can be expressed in this numerical representation as $r_1 = r_N = 0$, $z_1=z_2$ and $z_{N-1}= z_N$. We collect all the nodal values in the vector $\mbs{P} = \{r_1, z_1, r_2, z_2, \dots , r_N, z_N    \}$. If local inextensibility is required, we need to discretize the field of Lagrange multipliers
\[
\lambda(u) \approx \sum_{J=1}^M \hat{B}_J \lambda_J,
\]
where $\hat{B}_J$ are taken here to be the quadratic B-Splines obtained from the same knot vector used to generate the functions $B_I$. We have checked numerically the stability and convergence of this velocity/pressure pair in the present context.

Plugging these representations into the different models, generically written as $\mathcal{L}$, and equating to zero the derivatives of $\mathcal{L}$ with respect to velocity of the control points $\{\dot{r}_I(t), \dot{z}_I(t) \}$, possibly $\lambda_J$ and the global multipliers $\mbs{\Lambda}$ (or equivalently evaluating the variations in equation (\ref{pvp}) at the B-Spline basis functions), the Galerkin semi-discrete equations follow as
\begin{eqnarray}
\mathsfb{D}(\mbs{P})\dot{\mbs{P}} + \mathsfb{L}(\mbs{P})\mbs{\lambda} & = & \mbs{f}(\mbs{P}) \\
\mathsfb{L}^T(\mbs{P}) \dot{\mbs{P}}  & = & \mbs{0} \nonumber
\label{semidiscrete}
\end{eqnarray}
where $\mbs{\lambda}$ collects all the Lagrange multipliers for a
given model. 

Specifically, recalling equation (\ref{D_vol}), the entries for the column of the constraint matrix $\mathsfb{L}$ corresponding to the global volume constraint are
\[
\mathsfi{L}^{\rm vol}_{Ir} = -C^{\rm vol}[B_I,0], ~~~ \mathsfi{L}^{\rm vol}_{Iz} = -C^{\rm vol}[0, B_I].
\]
Similarly, recalling equation (\ref{D_area}), the column of  $\mathsfb{L}$ corresponding to the global area constraint is
\[
\mathsfi{L}^{\rm area}_{Ir} = -C^{\rm area}[B_I,0], ~~~ \mathsfi{L}^{\rm area}_{Iz} = -C^{\rm area}[0, B_I].
\]
The matrix entries of the local inextensibility constraints follow from Eq.~(\ref {inext_axi_lagr}):
\[
\mathsfi{L}^{\rm inext}_{Ir,J}  = - \int_0^1 \left(\frac{r'}{a^2}B_I' + \frac{1}{r} B_I\right) \hat{B}_J (2\pi ar){\rm d}u, ~~~ \mathsfi{L}^{\rm inext}_{Iz,J}  = - \int_0^1 \frac{z'}{a^2}  B_I'  \hat{B}_J (2\pi ar){\rm d}u.
\]
Recalling equation (\ref{G_rz}), the vector of nodal elastic forces is computed as
\[
f_{Ir} = G[B_I,0], ~~~ f_{Iz} = G[0, B_I].
\]
From equation (\ref{Wmem}), it follows that the only nonzero entries of $\mathsfb{D}^{\rm mem}$ are
\[
\mathsfi{D}^{\rm mem}_{Ir,Jr} = \int_0^1 \mu B_I B_J \frac{8 \pi a}{r} {\rm d}u.
\]

\subsection{The bulk dissipation matrix $\mathsfb{D}^{\rm bulk}$}

To discretize the bulk dissipation matrix, we follow a Galerkin Boundary Element method based on B-Splines. The traction jumps are discretized using the same B-Spline basis functions
\begin{equation}
\Delta \mbs{g}(u) \approx \sum_{I=1}^N B_I(u) \Delta \mbs{g}_I.
\label{dg}
\end{equation}
Recalling equations (\ref{vh}) and (\ref{dg}), multiplying equation (\ref{bir}) by a test function $B_I$ and integrating over the surface, we have
\[
\underbrace{\delta_{ij} \int_\Gamma  B_I  B_J {\rm d}S}_{\mathsfi{M}_{IiJj}} ~\dot{{P}}_{Jj} =\underbrace{- \frac{1}{\mu^{\rm bulk}}  \left(\int_\Gamma B_I(u) \int_\Gamma \mathsfi{G}_{ik}(\mbs{c}(u),\mbs{c}(v))   B_K(v) {\rm d}S(v) {\rm d}S(u)\right)}_{\mathsfb{A}_{IiKk}} \Delta{g}_{Kk},
\]
where summation over repeated indices is implied, and the lower-case indices $i$, $j$ and $k$ run over $r$ and $z$. Let us write the discretized dissipation potential for the bulk fluid. We introduce the vector
\[
\mbs{G}= \{\Delta g_{1r}, \Delta g_{1z},   \Delta g_{2r}, \Delta g_{2z}, \dots, \Delta g_{Nr}, \Delta g_{Nz} \}
\]
collecting the nodal traction jumps. It must be carefully noted that, unlike the vector of nodal forces $\mbs{f}$ and the other terms in the first Equation in (\ref{semidiscrete}), $\mbs{G}$ is not power-conjugate to $\dot{\mbs{P}}$. We have
\[
W^{\rm bulk} = \frac{1}{2} \int_\Gamma \mbs{V}\bcdot \Delta \mbs{g}~{\rm d}S = \frac{1}{2} \sum_I \sum_J \left(\int_\Gamma B_I B_J ~{\rm d}S\right)  \Delta \mbs{g}_I \bcdot \dot{\mbs{P}}_J = \frac{1}{2}  \mbs{G}^T \mathsfb{M} \dot{\mbs{P}},
\]
which, recalling that $\mbs{G} = \mathsfb{A}^{-1}\mathsfb{M}\dot{\mbs{P}}$ , becomes
\[ 
W^{\rm bulk} = \frac{1}{2} \dot{\mbs{P}}^T \mathsfb{M}\mathsfb{A}^{-1} \mathsfb{M} \dot{\mbs{P}},
\] 
and we finally identify $\mathsfb{D}^{\rm bulk} = \mathsfb{M}\mathsfb{A}^{-1} \mathsfb{M}$.

\subsection{Implementation details}

The computational strategy described above has been implemented within Matlab, except for the boundary integral method, which has been coded in \texttt{C++} and is built on top of
the \texttt{deal.II} library \citep{deal-II} with B-Spline support based on the \texttt{GSL} library \citep{gsl}. A Matlab interface of this code has been set up. This system of differential-algebraic equations in (\ref{semidiscrete}) is integrated in time with the specialized solvers \texttt{ode15s} and \texttt{ode23t}.

Special care must be taken in performing the integration of the
diagonal terms of the $\mathsfb{A}$ matrix, since the axisymmetric
Stokeslet~(\ref{axisymmetric-stokes-kernel}) presents a logarithmic
singularity when its arguments are close to each other. We use integration formulas based on weighted Gauss quadrature
rules. 

\subsection{Reparameterization of the B-Spline curve}

\begin{figure}
\begin{center}
\includegraphics[width=6.5cm]{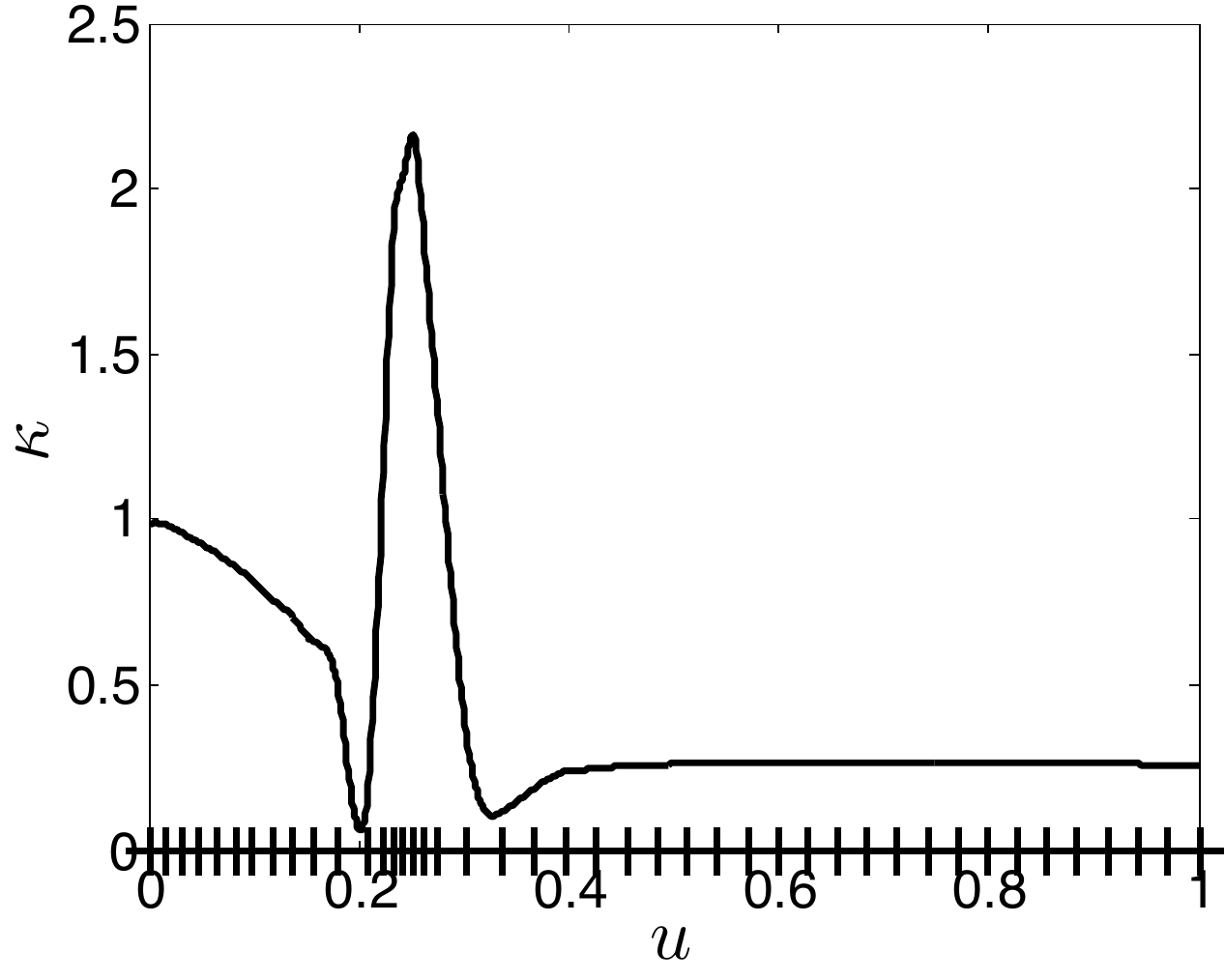}
\includegraphics[width=6.5cm]{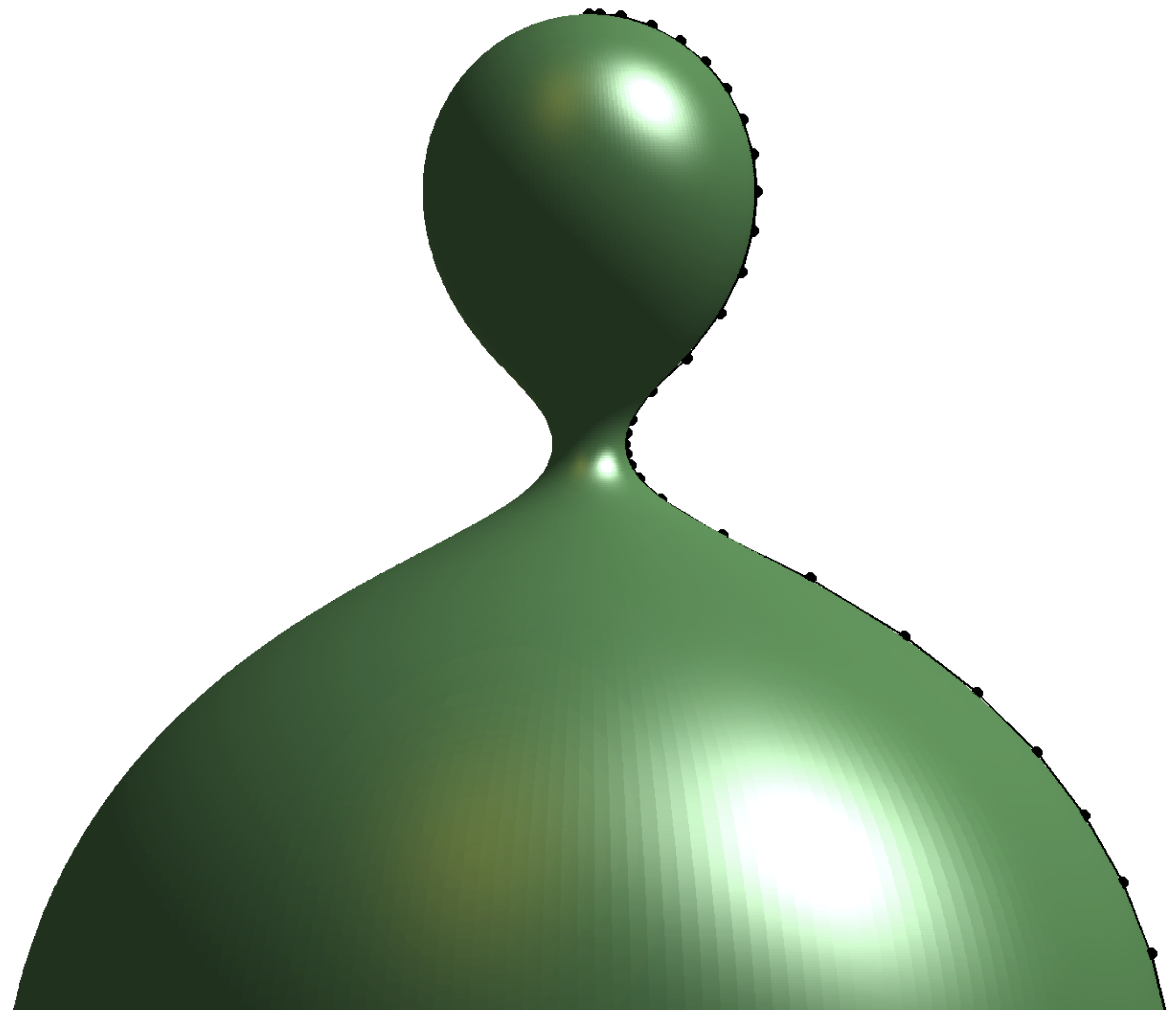}
\caption{Non-uniform knot span adapted to the curvature of the generating curve (left), and control polygon for a given vesicle shape (right). The the parameterization of $\mbs{c}(u)$ is nearly arc-length. The clustering in the knots of the B-Spline results in a clustering of the control points in highly curved regions.}
\label{reparam}
\end{center}
\end{figure}

As the system in equation (\ref{semidiscrete}) is integrated in time, the control points of the generating curve collected in $\mbs{P}$ may cluster in some areas, over-resolving some parts of the domain and leaving others with poor resolution, irrespective of the geometric features of the generating curve. Furthermore, as illustrated in the examples later, the highly curved parts of the generating curve may move significantly in the parametric domain $[0,1]$. For these reasons, it is important to re-parameterize the generating curve periodically in the simulations. 

Here, we have devised a heuristic method that builds a non-uniform knot-span for the B-Spline functions, which clusters knots, hence basis functions, where the curvature of the generating curve is high. At the re-parameterization instants, the new control points relative to the new set of basis functions are obtained by a least square fit to the previous description of the generating curve that results in a nearly arc-length parameterization. See Figure \ref{reparam} for an illustration. This method allows us to  reduce substantially the number of control points as compared to resolving the fine geometrical features with uniform refinement, and still obtain very accurate solutions.

\section{Results}
\label{results}

\begin{figure}
\begin{center}

\includegraphics[height=6.9cm]{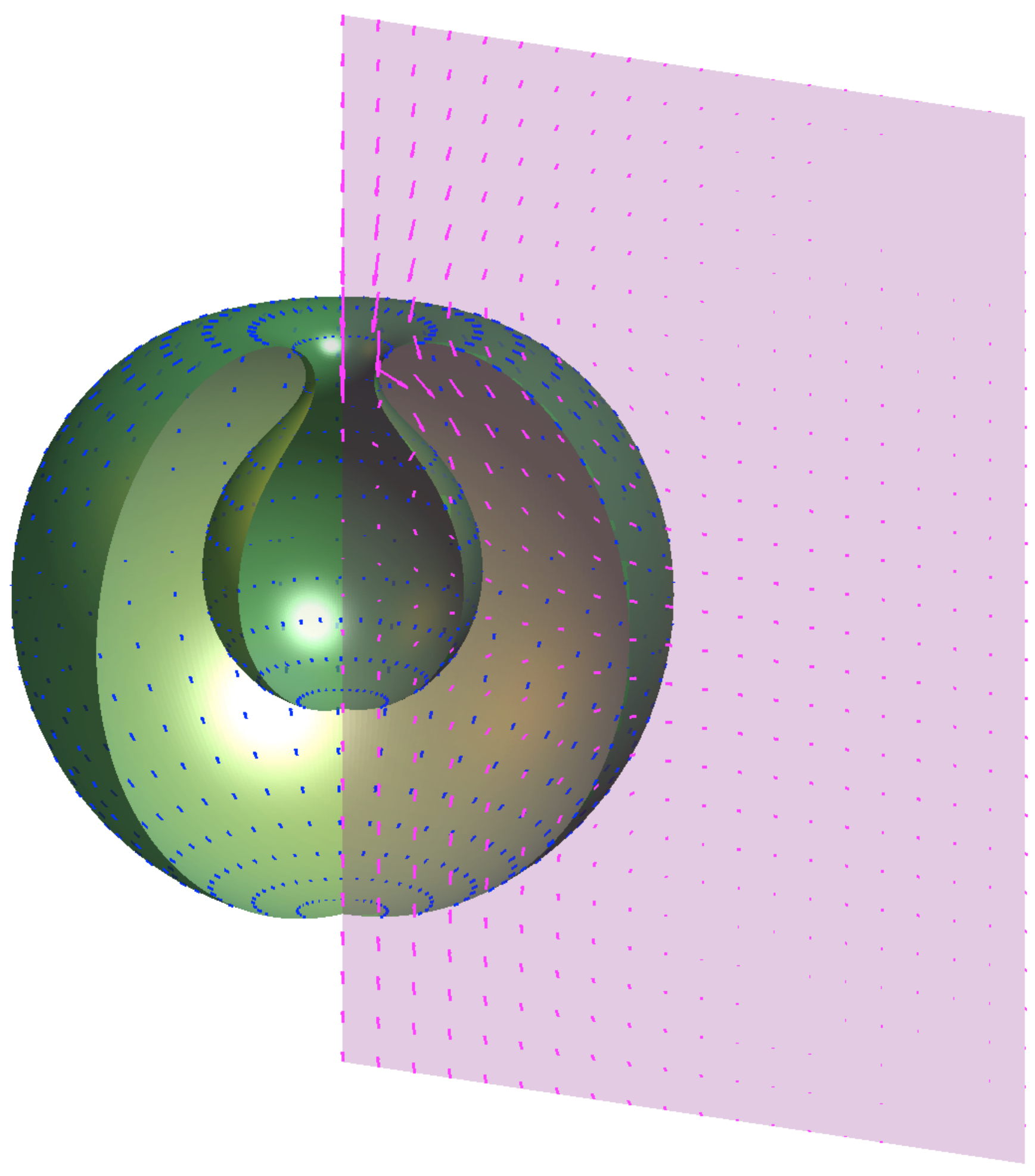}
\includegraphics[height=6.9cm]{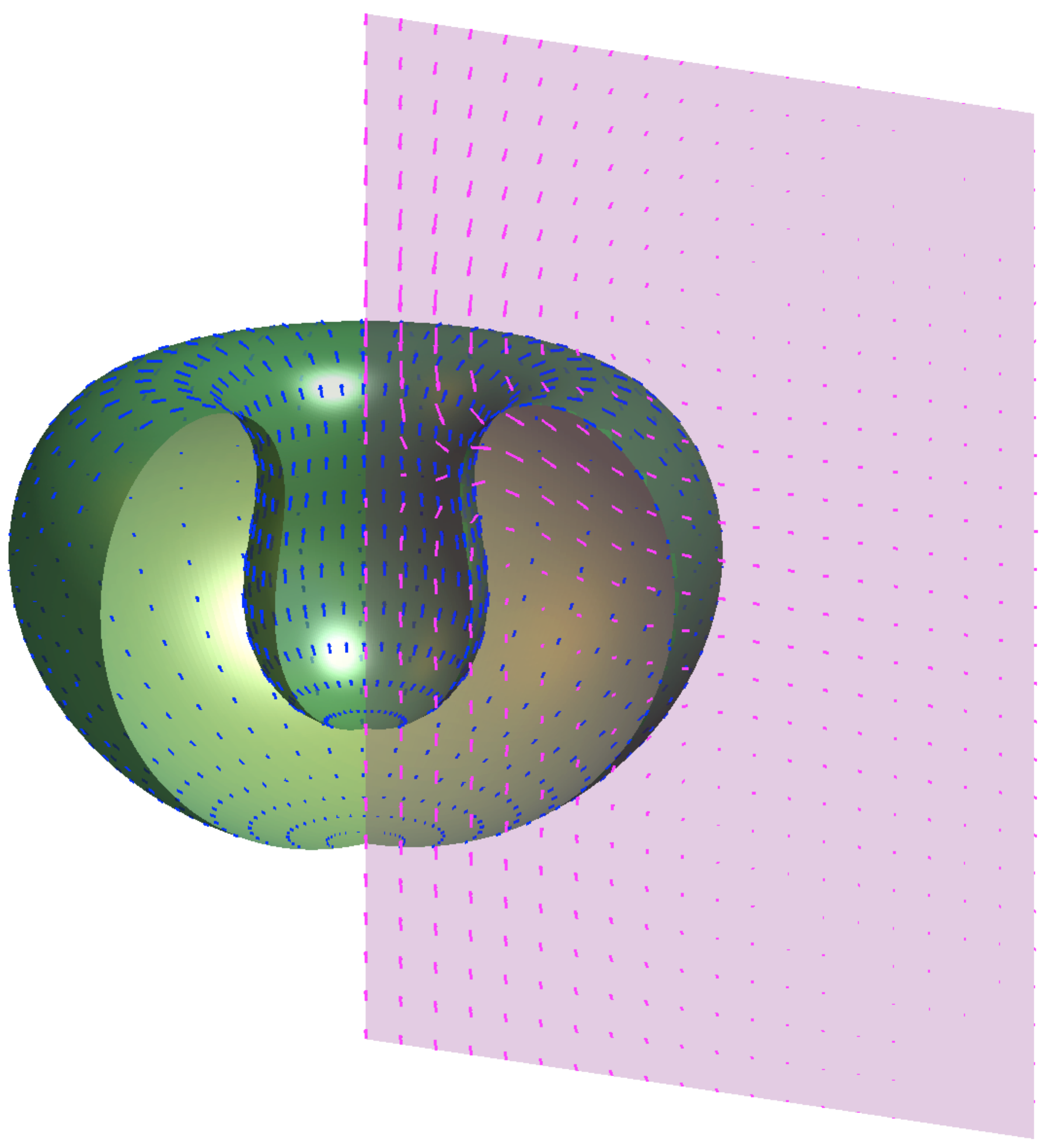}
\includegraphics[height=6.9cm]{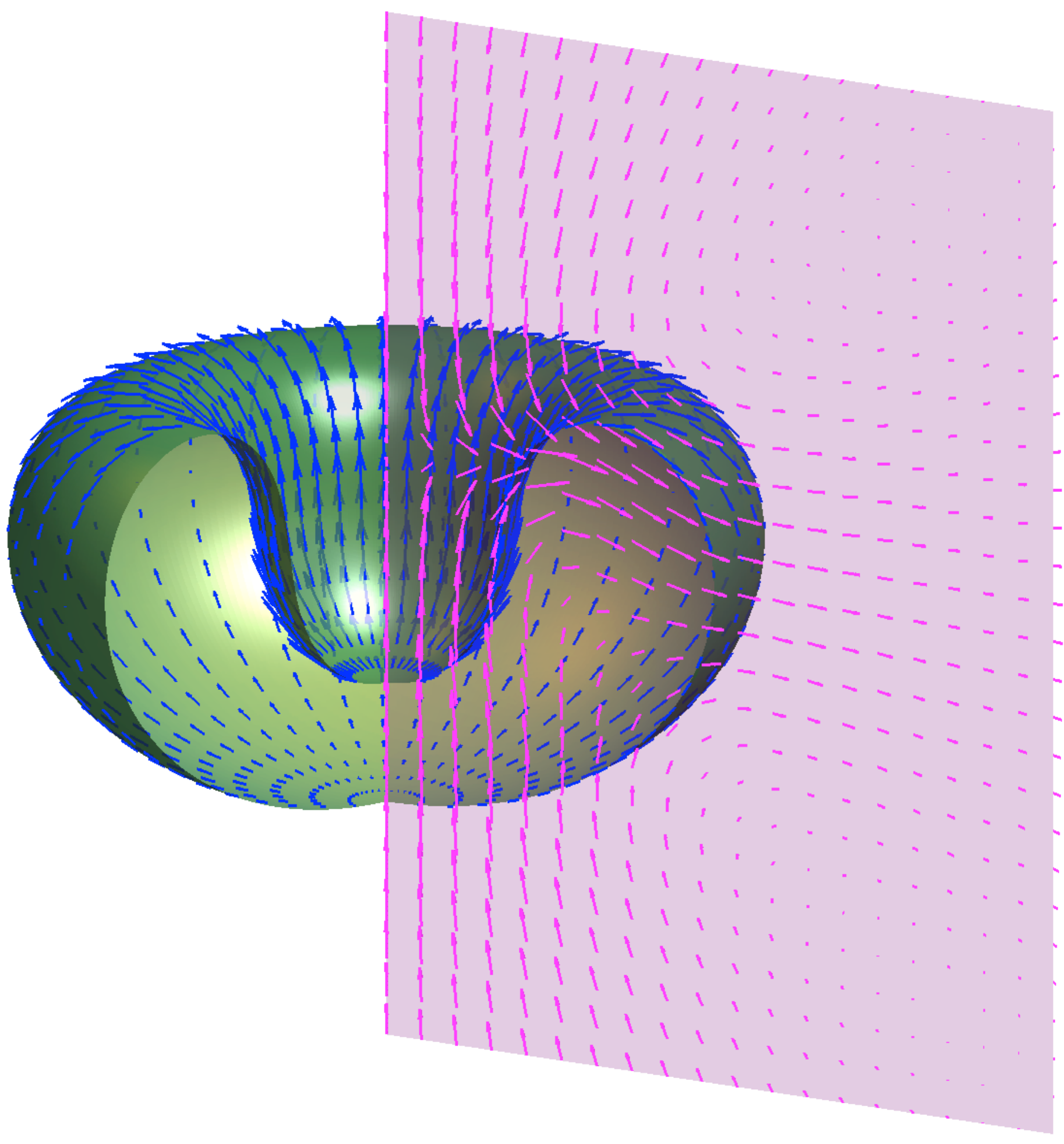}
\includegraphics[height=6.9cm]{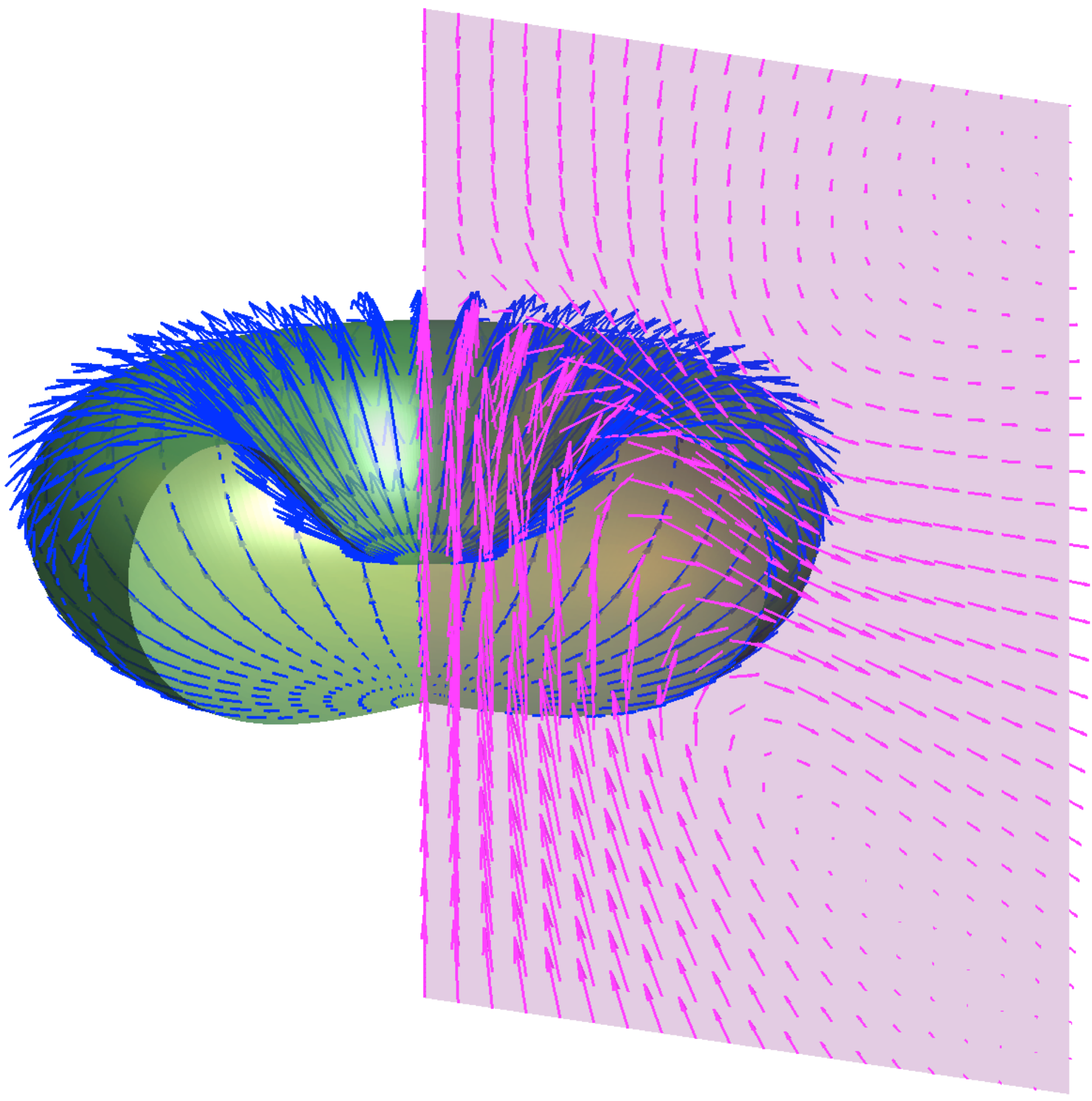}
\includegraphics[height=6.8cm]{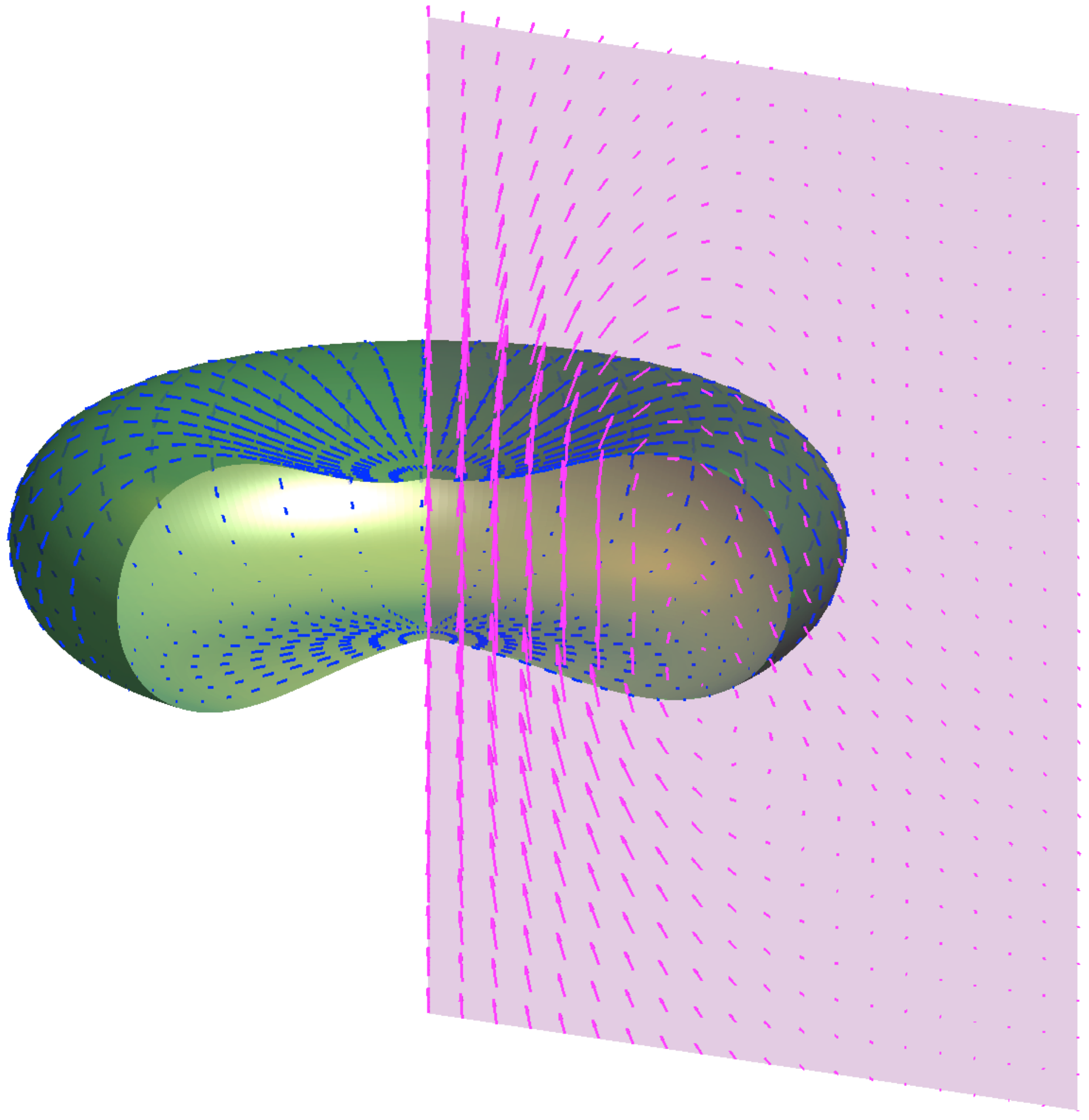}
\includegraphics[height=6.8cm]{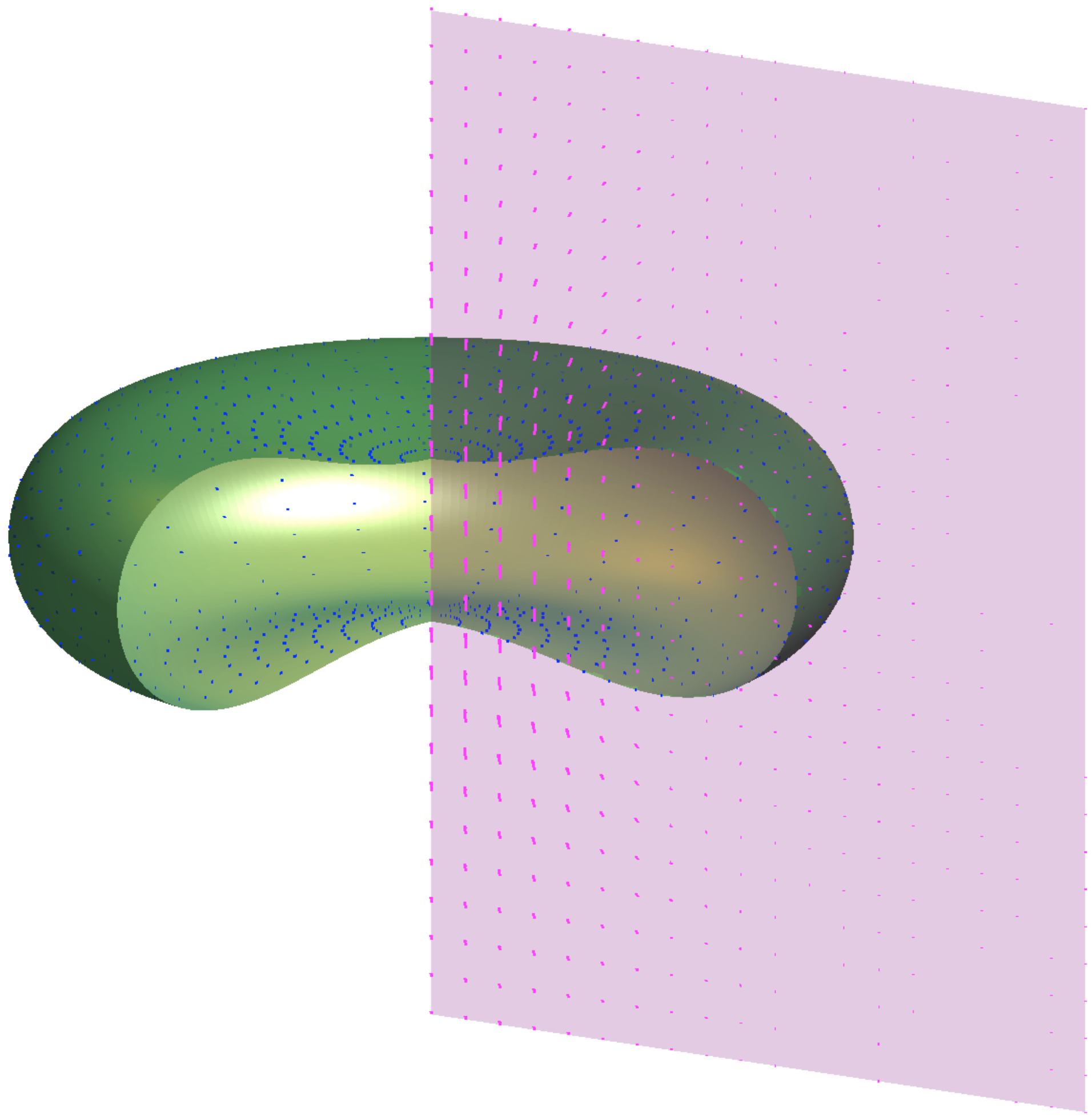}

\caption{Snapshots of a typical simulation. Sector of the axisymmetric shape (green), tangential velocity field on the membrane representing the flow of amphiphiles as they rearrange to accommodate the shape evolution (blue), and velocity field of the ambient fluid post-processed from the boundary integral solution (magenta).}
\label{typical}
\end{center}
\end{figure}

Snapshots from a typical simulation are shown in figure \ref{typical}. The tangential velocity on the membrane is represented to visually illustrate the inextensible flow of amphiphiles on the membrane, whose tangential shear causes the membrane dissipation. The ambient fluid velocity field has also been represented, as post-processed from the boundary integral numerical method. This bulk velocity field is incompressible, its tangential projection on the membrane is inextensible, and its normal component on the membrane illustrates the rate of change of the shape.

We next report two different sets of relaxation dynamics simulations. On the one hand, the qualitative differences of the flows engendered by the different dissipative mechanisms are analyzed. On the other hand, we consider the competition between the bulk and the membrane dissipative mechanisms in $\mathcal{L}^{\rm full}$ as a function of the system size, and obtain a composite power law for the relaxation time as a function of size.

\subsection{Comparison of the dynamics for different dissipation mechanisms}

We consider in this section a gallery of relaxation simulations, and compare qualitatively the dynamics obtained with each individual dissipation mechanism, i.e. the membrane dissipation in model $\mathcal{L}^A$, the bulk dissipation in model  $\mathcal{L}^B$, and the $L_2$ dissipation in  model $\mathcal{L}^C$. In the examples in this section, the time scale is arbitrary and has been re-scaled so that, for a given test case, all the models reach a given energy threshold at the same instant.  

\subsubsection{Relaxation of a discocyte vesicle}

\begin{figure}
\begin{center}
\includegraphics[width=13.cm]{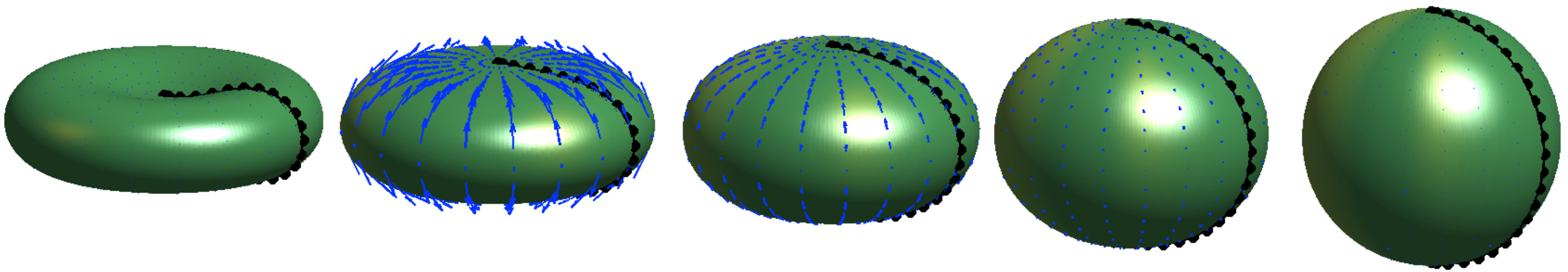}
\includegraphics[width=5.cm]{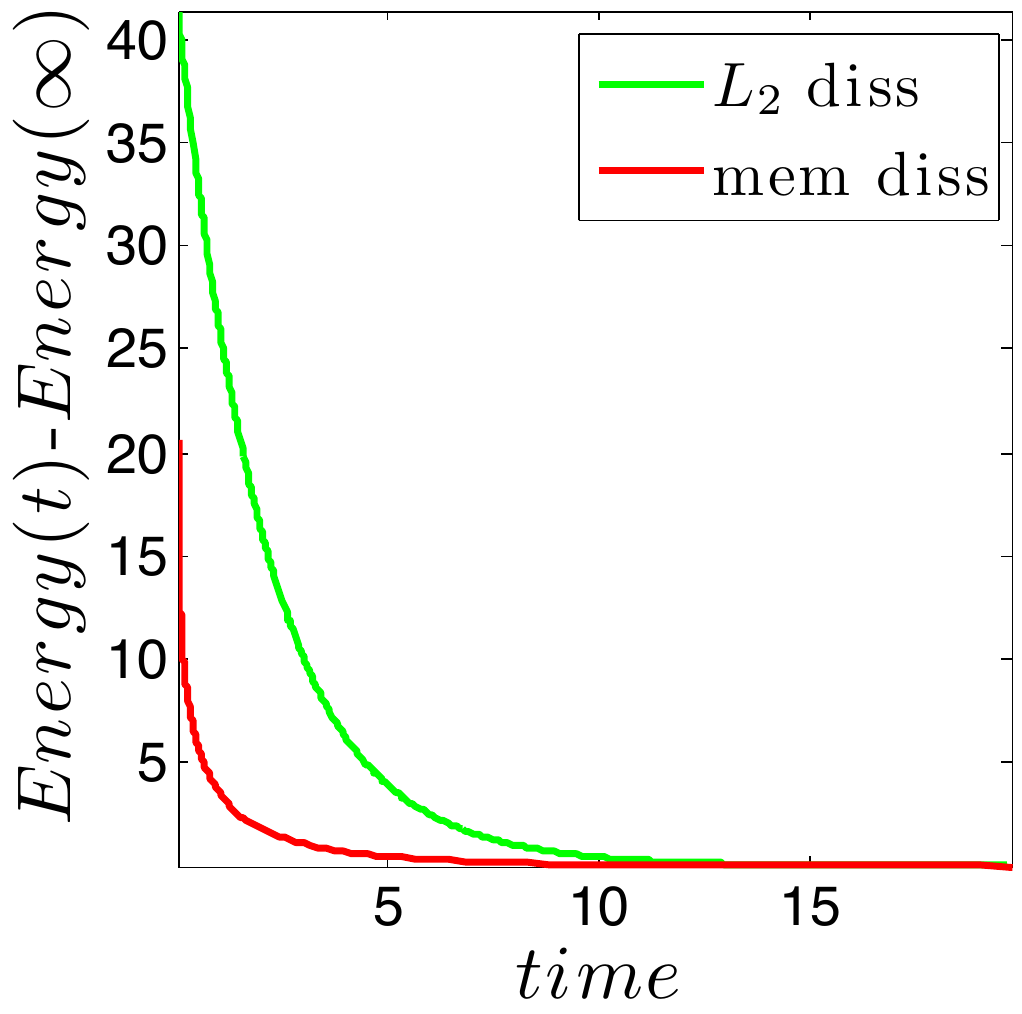}
\includegraphics[width=5.cm]{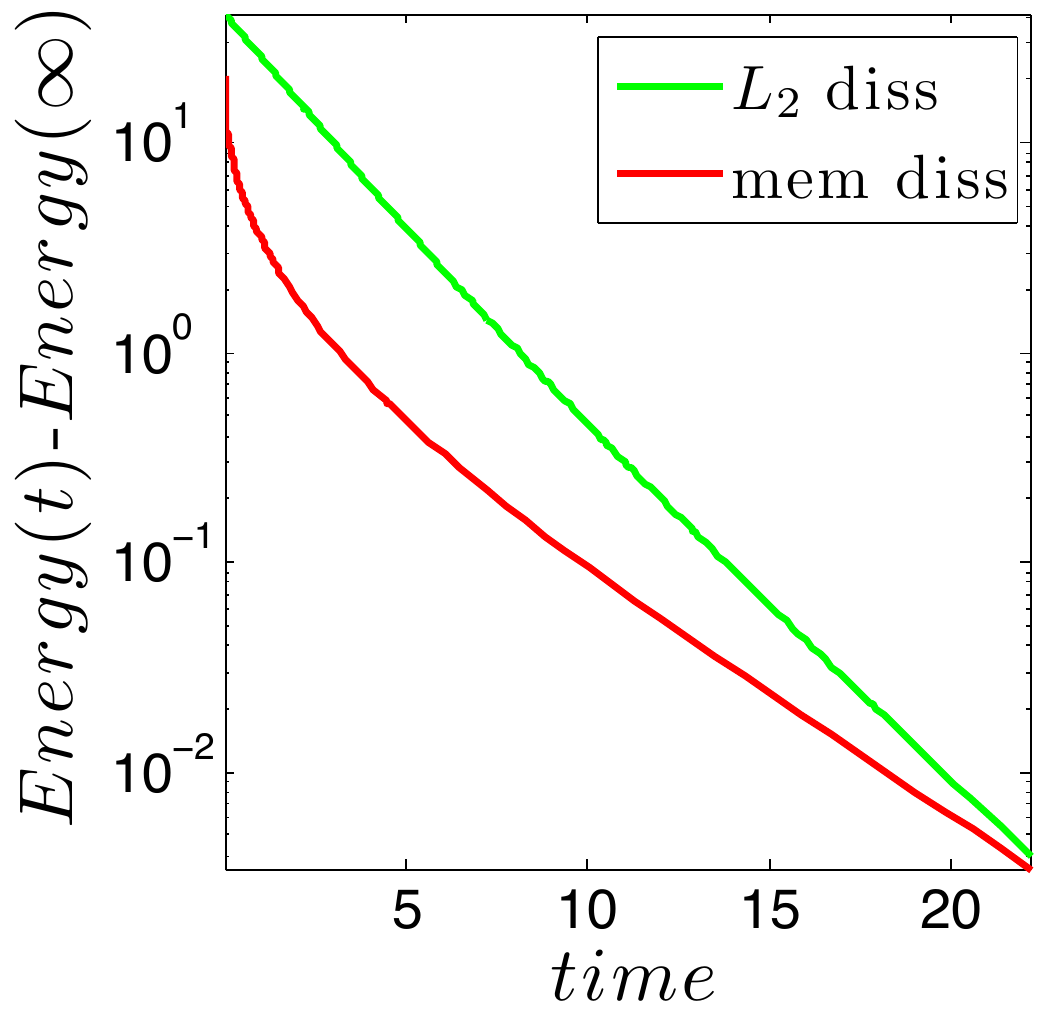}
\caption{Discocyte vesicle evolving towards the sphere after the volume constraint is released. Snapshots of the evolution (top), and elastic energy of the membrane as a function of time in linear and semi-logarithmic scales, for models $\mathcal{L}^A$ and $\mathcal{L}^C$ (bottom).}
\label{disco}
\end{center}
\end{figure}

We consider first the relaxation of a discocyte vesicle, in equilibrium for zero spontaneous curvature and a reduced volume of 
\[
v = 3 \sqrt{4\pi} \frac{V}{S^{3/2}} = 0.65,
\]
where $V$ is the volume of the vesicle and $S$ its surface area  \citep{Handbook-Biological-Physics}. The reduced volume is 1 for a sphere, and smaller for any other shape. In this somewhat artificial example, the volume constraint is suddenly released, and the system evolves from the discocyte configuration towards the sphere adopting oblate configurations along the way. For this example it does not make sense to consider the bulk fluid viscosity, i.e. functional $\mathcal{L}^B$, since the enclosed volume cannot be incompressible. Figure \ref{disco} illustrates the relaxation dynamics, as well as the time evolution of the elastic energy. This evolution is presented also in semi-logarithmic scale, which better shows the qualitative behavior in that it is easy to visually detect deviations from the exponential relaxation of a simple linear system characterized by a single time-scale. In this example, model $\mathcal{L}^C$ without volume constraint behaves nearly exponentially, as shown by the linear response of this model in the semi-logarithmic plot. The sharp deviation from the exponential behavior of model $\mathcal{L}^A$ without volume constraint is evident from the plot, particularly at early stages. As a matter of fact, the energy evolution for the membrane dissipation model seems to have a vertical asymptote at the origin.

\begin{figure}
\begin{center}
\includegraphics[width=7.cm]{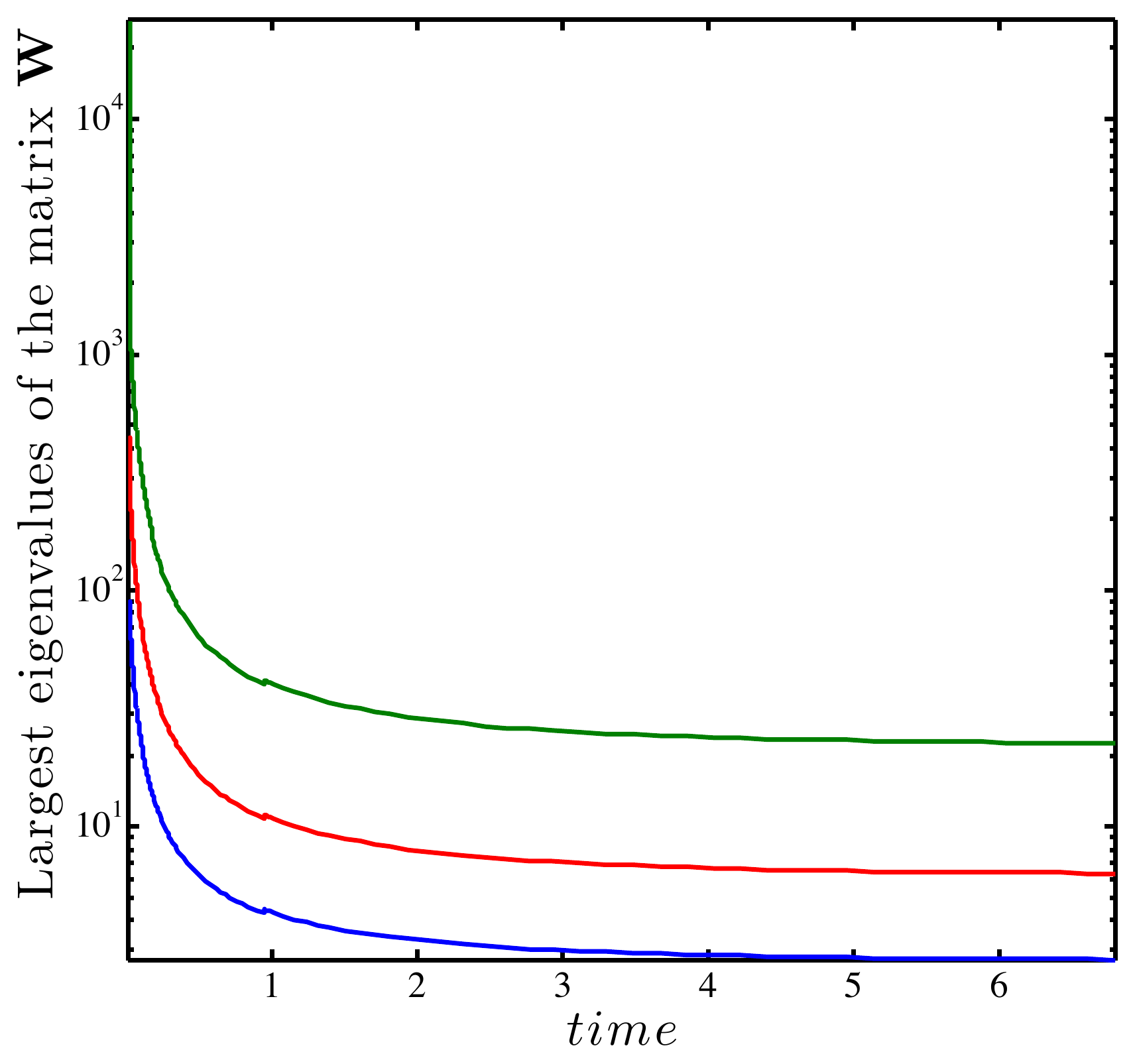}
\caption{Time evolution of the first three eigenvalues of the matrix $\mathsfb{W}$ for the relaxing discocyte vesicle.}
\label{eigs}
\end{center}
\end{figure}

To investigate the origin of this behavior, we eliminate the Lagrange multipliers from equation (\ref{semidiscrete}) by isolating $\dot{\mbs{P}}$ in the first equation, substituting in the second to isolate $\mbs{\lambda}$, and substituting back in the first equation, which results in
\[
\dot{\mbs{P}} = \underbrace{\left(\mathsfb{D}^{-1} -  \mathsfb{D}^{-1}  \mathsfb{L} (\mathsfb{L}^T {\mathsfb{D}^{-1}} \mathsfb{L})^{-1} \mathsfb{L}^T   \mathsfb{D}^{-1}     \right)}_{\mathsfb{W}(\mbs{P})} \mbs{f}(\mbs{P}).
\]
The largest eigenvalues of the matrix $\mathsfb{W}\left(\mbs{P}(t)\right)$ characterize the effective inverse viscosity of the dominant modes of the system, once the constraints have been eliminated. While for the model based on the $L_2$ dissipation, $\mathcal{L}^C$, the eigenvalues remain nearly constant during the time evolution, figure \ref{eigs} shows that the maximum eigenvalue for the model based on the membrane dissipation, $\mathcal{L}^A$, changes by three orders of magnitude. This shows the strong geometry dependence of the membrane dissipation, which seems to nearly vanish for discocytes, and then rapidly increases in the transition to oblate morphologies. We revisit later this behavior in the example of a relaxing stomatocyte.

\subsubsection{Relaxation of a pearled vesicle}

\begin{figure}
\begin{center}
\begin{minipage}[b]{7.4cm}
\includegraphics[width=7.4cm]{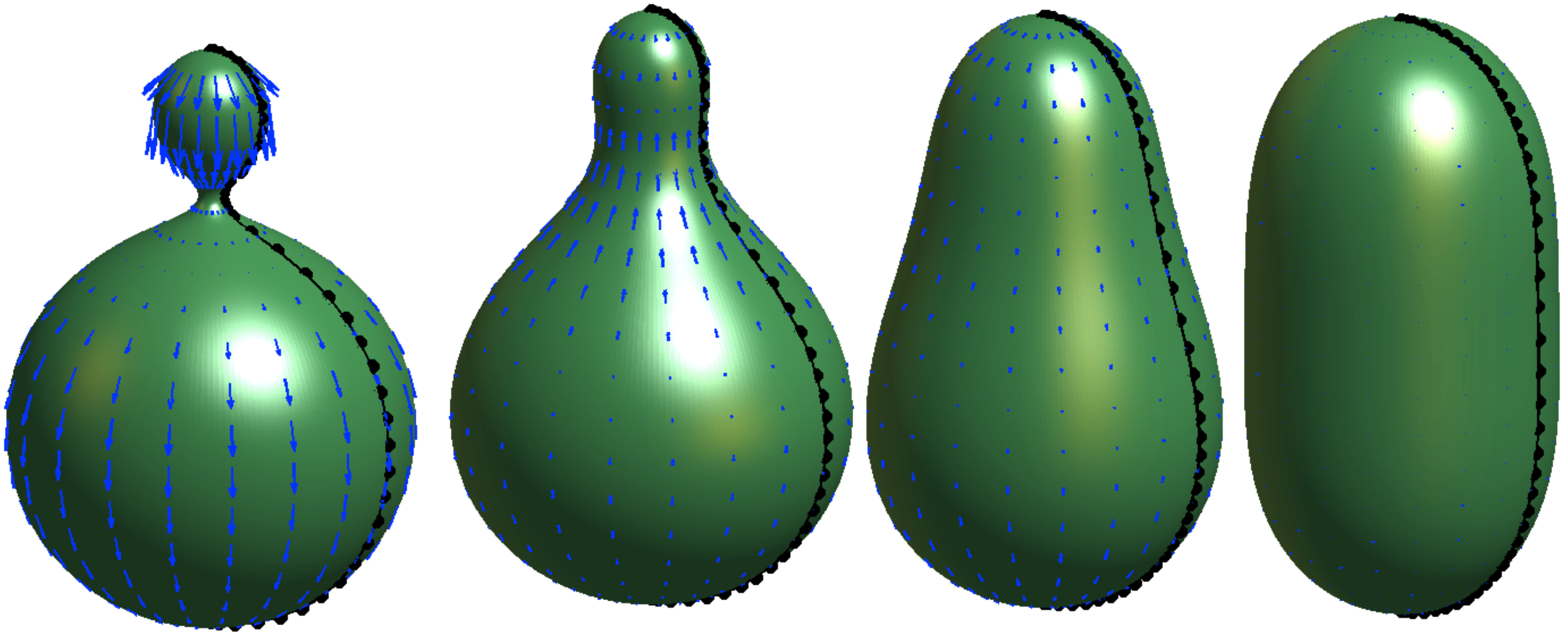}
\\
\\
\\
\\
\end{minipage}
\begin{minipage}[t]{6cm}
\includegraphics[width=6.cm]{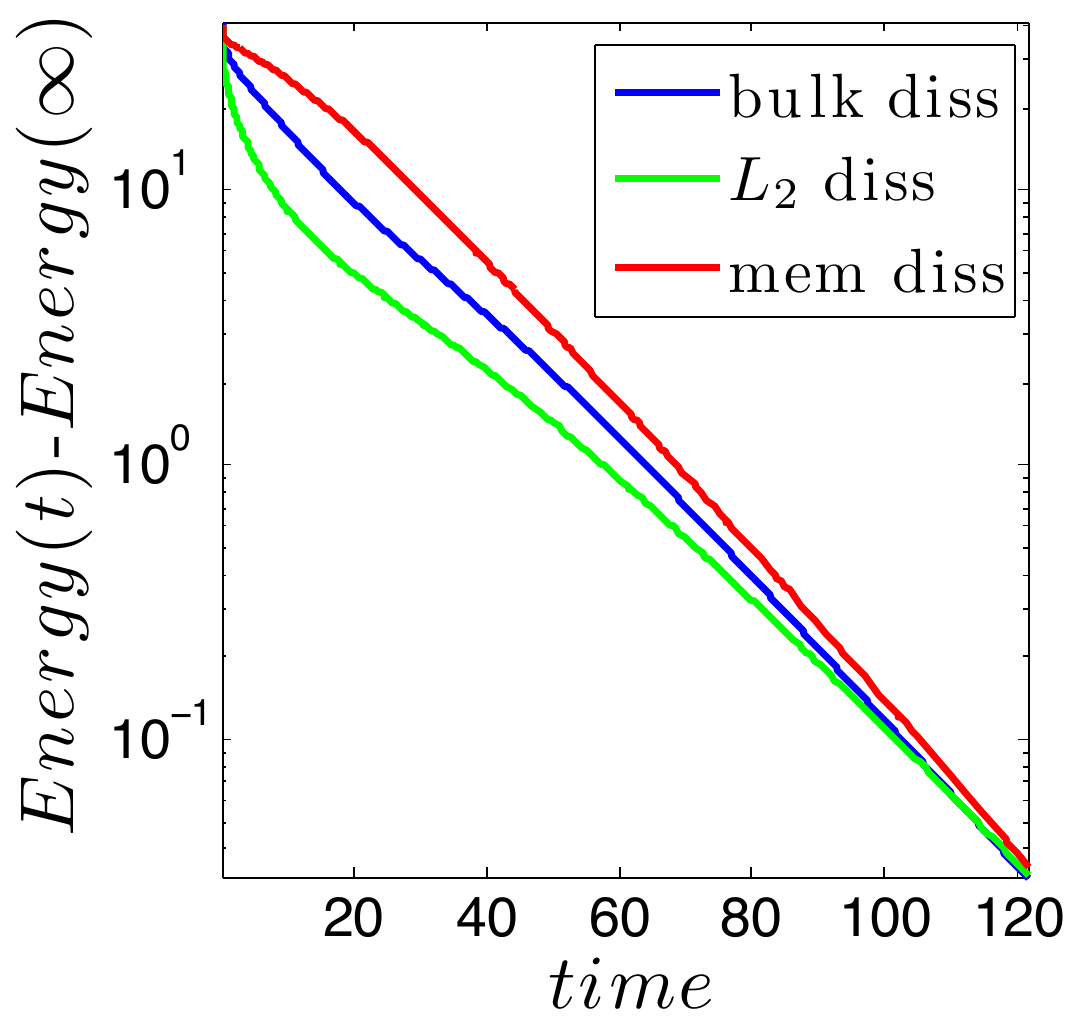} 
\end{minipage}
\caption{Pearled vesicle evolving towards a prolate configuration, modeling the large scale relaxation after a fusion event between two dissimilar vesicles. Snapshots of the evolution (left), and elastic energy of the membrane as a function of time in semi-logarithmic scale, for models $\mathcal{L}^A$, $\mathcal{L}^B$, and $\mathcal{L}^C$ (right).}
\label{examples_pearled}
\end{center}
\end{figure}

We now consider the relaxation of a pearled vesicle with the enclosed volume constraint in place. The initial configuration is in equilibrium for a non-zero value of the spontaneous curvature. In the dynamical simulation, $C_0$ is set to zero, and then the system evolves towards a prolate configuration, see Figure \ref{examples_pearled} (left). The experiment can be viewed as the large scale shape evolution after two dissimilar vesicles have fused and the fusion pore has expanded to some extent. The energy relaxation with the three models, $\mathcal{L}^A$, $\mathcal{L}^B$, and $\mathcal{L}^C$, is shown in Figure \ref{examples_pearled} (right). For this example, in contrast with the example above, the $L_2$ evolution is the one which deviates most from the exponential relaxation, with an initial very fast regime, for which the system can significantly lower its energy without much dissipation, then a slowdown and finally an exponential tail. The membrane and the bulk dissipation evolutions deviate slightly from the exponential relaxation, the former with a slower evolution and the latter with a faster evolution at early times.

\subsubsection{Relaxation of a tethered vesicle}

\begin{figure}
\begin{center}
\includegraphics[width=11.cm]{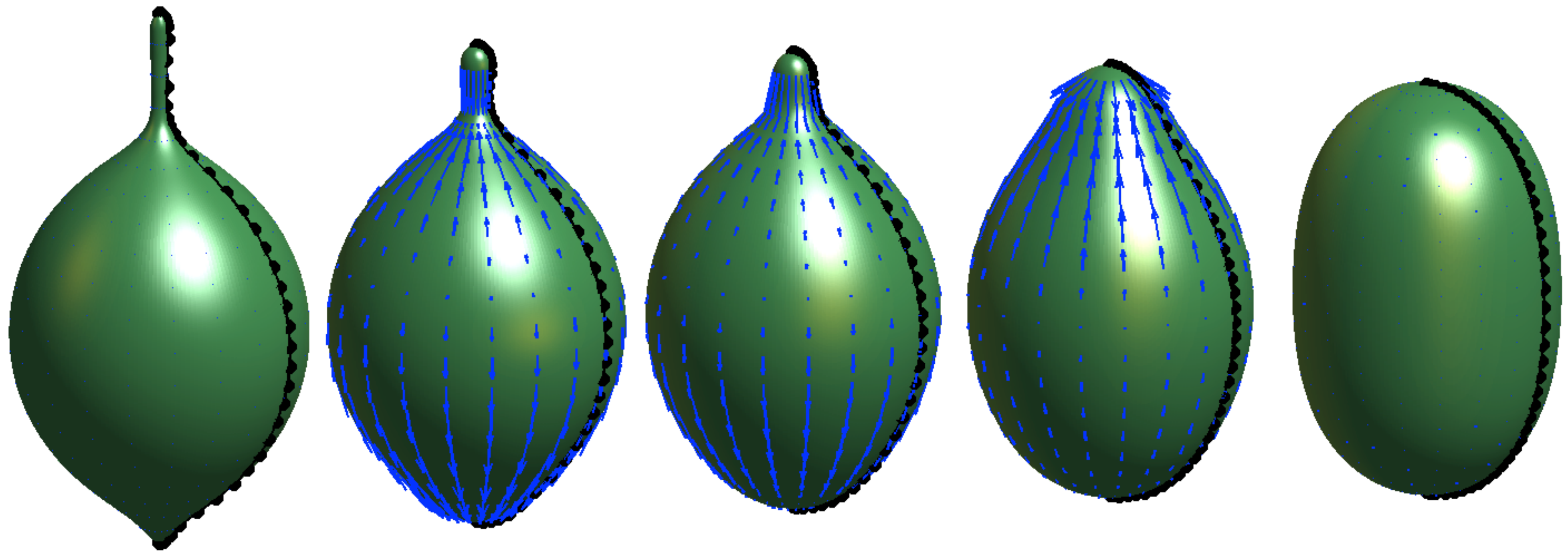}
\includegraphics[width=5.cm]{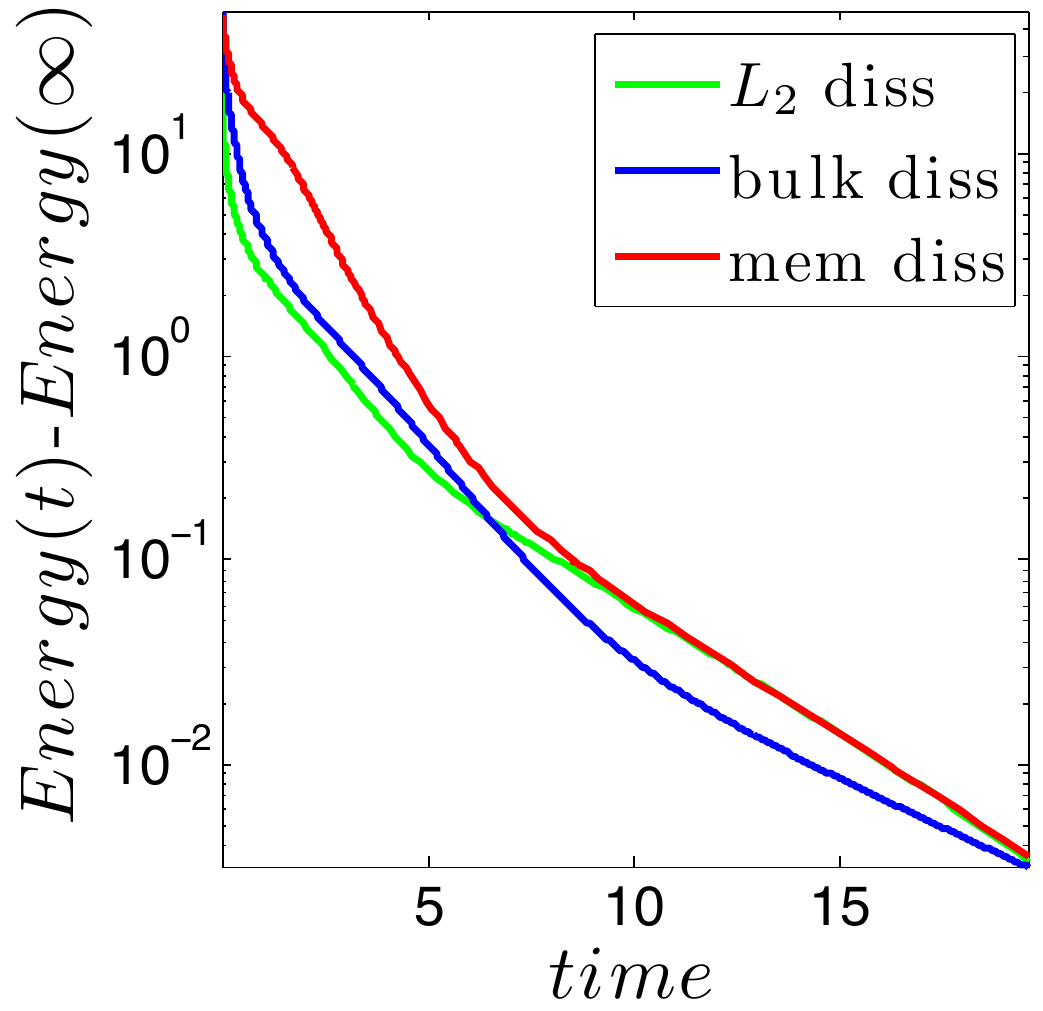}
\caption{Tethered vesicle evolving towards a prolate configuration, after the force keeping the vesicle stretched is released. Snapshots of the evolution (top), and elastic energy of the membrane as a function of time in semi-logarithmic scale, for models $\mathcal{L}^A$, $\mathcal{L}^B$, and $\mathcal{L}^C$ (bottom).}
\label{example_tether}
\end{center}
\end{figure}

We now study the relaxation dynamics of a tethered vesicle, prepared computationally following the experimental setup in \cite{phillips-tethers} in which a vesicle is pulled by two beads controlled by optical tweezers. A breaking of the symmetry was observed experimentally, and also in our simulations, as in one of the attachments the vesicle produces a tether to relieve the elastic energy while maintaining the volume and area constraints. Here, we report the relaxation as the forces stretching the vesicle are released, see Figure \ref{example_tether}. Not suprisingly, the features of the relaxation dynamics are similar to those in the previous example. However, in this case, two exponential regimes (linear in the semi-logarithmic scale) can be observed, corresponding to the relatively fast event of the tether retraction, and to the slower event of the global convergence to the prolate configuration.

\subsubsection{Relaxation of a stomatocyte vesicle}

\begin{figure}
\begin{center}
\includegraphics[width=13.cm]{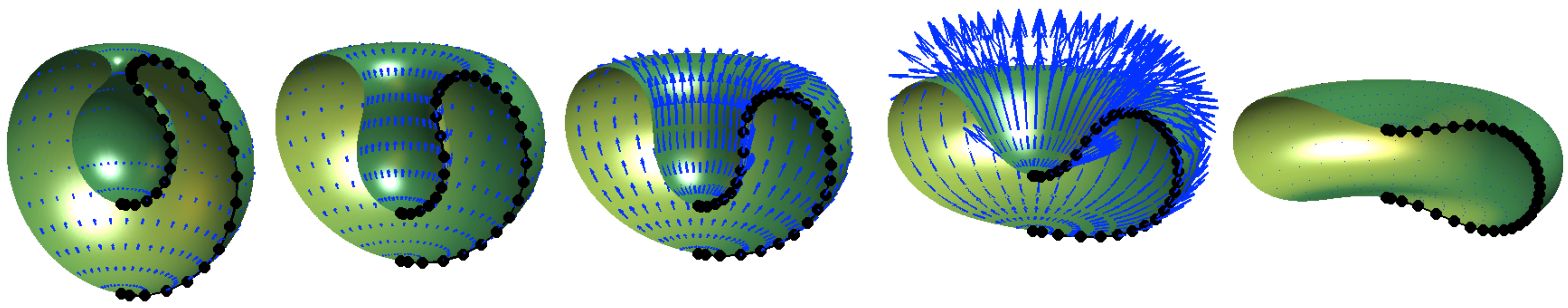}
\includegraphics[width=5.cm]{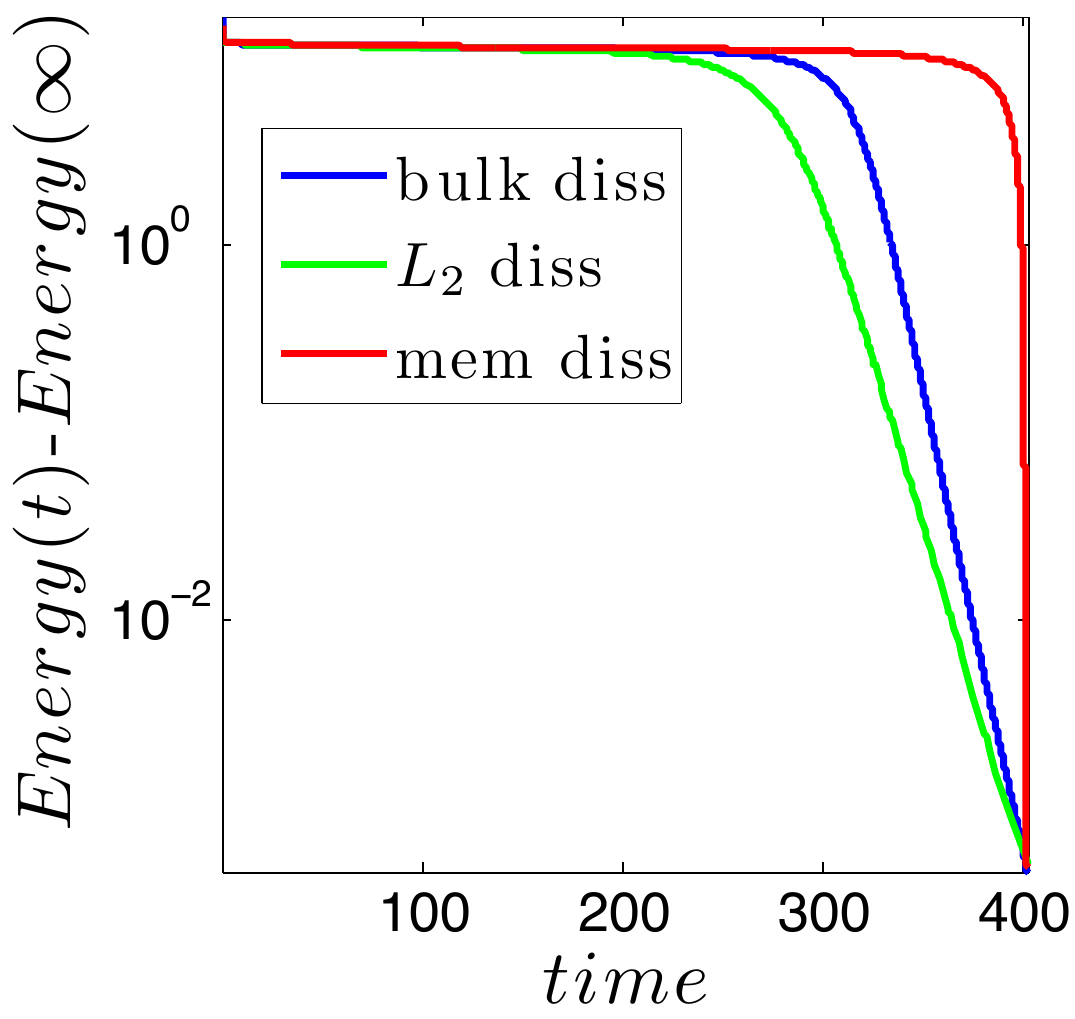}
\caption{Stomatocyte vesicle evolving towards a discocyte configuration, modeling the relaxation after the fusion of two dissimilar vesicles. Snapshots of the evolution (top), and elastic energy of the membrane as a function of time in semi-logarithmic scale, for models $\mathcal{L}^A$, $\mathcal{L}^B$, and $\mathcal{L}^C$ (bottom).}
\label{example_stomato}
\end{center}
\end{figure}

We consider here the relaxation of a stomatocyte vesicle, prepared by minimizing the elastic energy with a non-zero spontaneous curvature. The experiment can be viewed as the large scale shape evolution after two dissimilar vesicles, one being inside of the other, have fused and the fusion pore has expanded to some extent (see Figure \ref{example_stomato}). This last example exhibits a dramatic dependence of the relaxation dynamics on the dissipative mechanism. Generically, we observe a very slow dynamics at early stages, where the elastic forces are nearly perpendicular to the constraint manifold. An exponential relaxation is observed for the three models. Once the neck of the stomatocyte is erased (third snapshot in  Figure \ref{example_stomato}, much faster dynamics are observed, with well-defined exponential tails characterized by much larger time constants in both the bulk dissipation and the $L_2$ dissipation models. In sharp contrast, the membrane dissipation model exhibits an even faster, non exponential relaxation at later stages, and abruptly relaxes to the equilibrium discocyte morphology. As a matter of fact, there seems to be a vertical asymptote in the energy vs time relation, a similar behavior to that observed in the first example, also near a discocyte configuration. This vanishing effective membrane viscosity for such configurations is worth analytical investigation, and highlights the strong geometry dependence of the membrane dissipation.

\subsection{Competition between membrane and bulk viscosity}

As a final more realistic example, we consider the relaxation dynamics when both the membrane and the bulk dissipative mechanisms are in place. In this case, a natural length-scale arises as the ratio between the membrane two-dimensional viscosity and the bulk viscosity \citep{arroyo-desimone-2009}
\begin{equation}
\ell = \frac{\mu}{\mu^{\rm bulk}}.
\label{visc_l}
\end{equation}
This length-scale arises also in the Saffman-Delbruck theory for the diffusion of membrane inclusions \citep{Saffman08011975}, see also \citep{Stone98}, and separates the behavior of large vesicles, for which the bulk viscosity dominates over the membrane viscosity, and that of small vesicles, whose dynamics are mostly dictated by the membrane dissipation. On the basis of a simple analytical model for budding, it has been estimated \citep{arroyo-desimone-2009} that for large vesicles, the relaxation time scales as
\begin{equation}
t_r \sim \frac{\mu^{\rm bulk} R_0^3}{\kappa},
\label{t_1}
\end{equation}
while for small vesicles, we have
\begin{equation}
t_r \sim \frac{\mu R_0^2}{2 \kappa}.
\label{t_2}
\end{equation}
The cross-over between the two power-laws are for vesicles sizes of $R_0 \sim \ell/2$, where $R_0$ is the characteristic size of the vesicle, here the radius of the half-sphere whose surface area is that of the vesicle under consideration.

\begin{figure}
\begin{center}
\includegraphics[width=9cm]{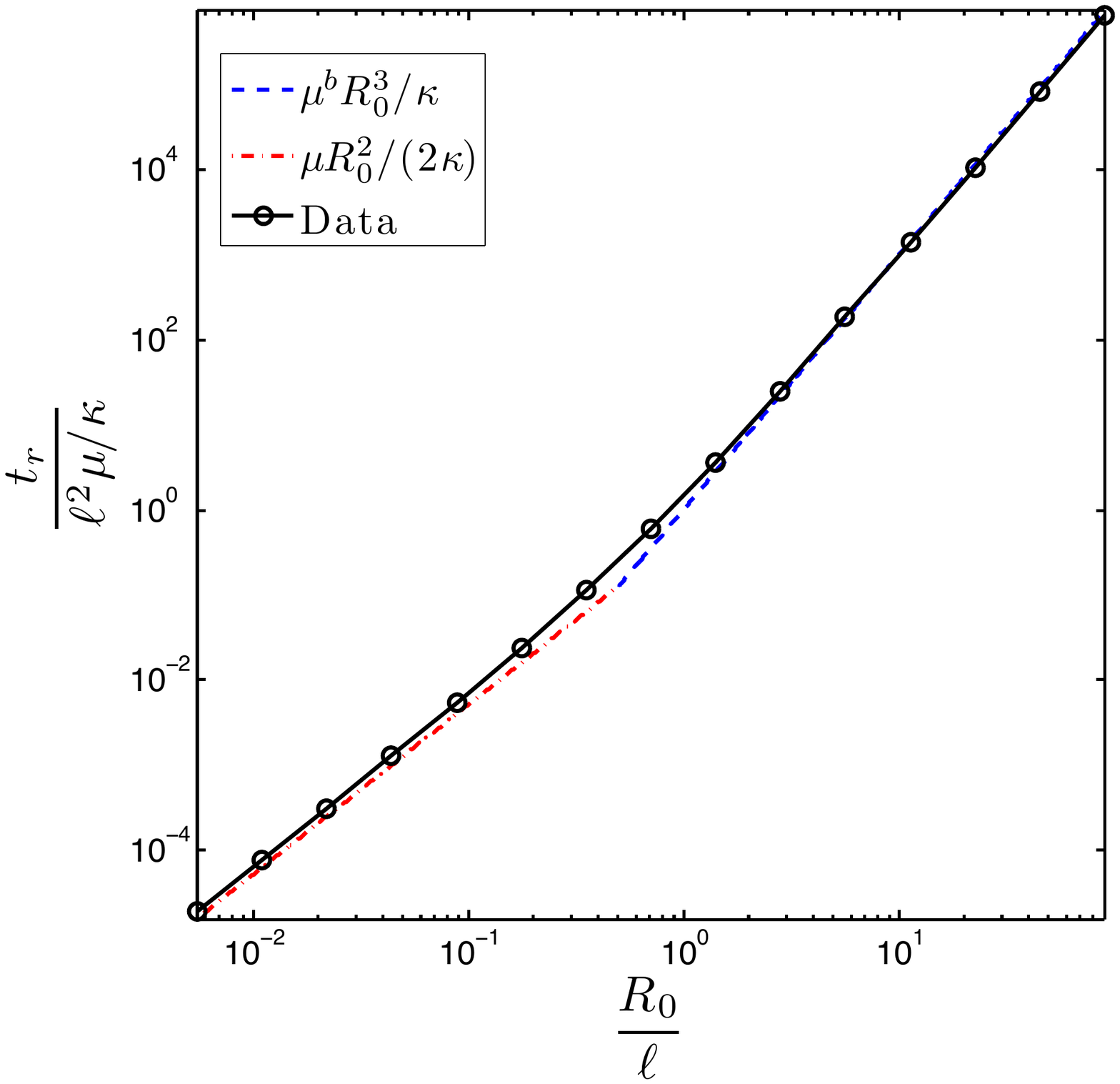}

\vspace{-4.5cm}
\hspace{5.4cm}
\includegraphics[height=2.9cm]{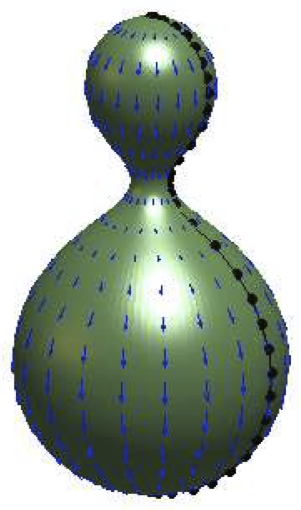}
\includegraphics[height=2.9cm]{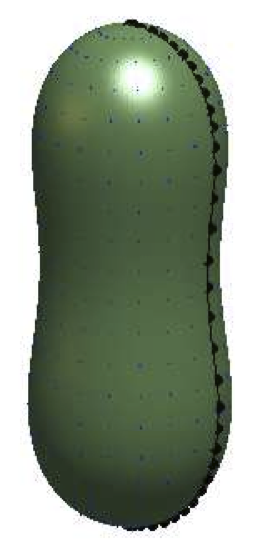}
\vspace{1.8cm}

\caption{Composite power law for the relaxation time as a function of vesicle size for a fusion event, with a regime dominated by bulk viscosity (large vesicles) and another dominated by membrane viscosity.}
\label{relaxation}
\end{center}
\end{figure}

We test these analytical estimates with a simulation resolving the actual shape relaxation coupled with the membrane and the bulk dissipations, as described by the functional $\mathcal{L}^{\rm full}$. We perform relaxation experiments of pearled vesicles of different sizes converging to an oblate configuration, and for each simulation, we record the relaxation time, defined as the time it takes for the system to relax up to 99.25\% of its excess elastic energy. The results are shown in a logarithmic scale in Figure \ref{relaxation}, together with the predictions of Eqs.~(\ref{t_1}-\ref{t_2}). The data of the simulations show a remarkable agreement with the composite power-law predicted by the analytical estimates, not only qualitatively, but also quantitatively in this example, even though the budding model in \cite{arroyo-desimone-2009} differs significantly from the relaxation simulations considered here.

\section{Conclusions}
\label{conclusions}

We have presented, to the best of our knowledge, the first calculations for the dynamics of fluid membranes accounting for the coupling between the morphological dynamics for finite shape changes and the membrane two-dimensional curved flow. To this end, we have particularized the theory in \cite{arroyo-desimone-2009} to the axisymmetric case, yielding a simple, workable model, which is the basis to the numerical simulations. The reported simulations show that the relaxation dynamics driven by curvature elasticity with the usual constraints are qualitatively different for the different dissipation mechanisms (membrane viscosity, bulk viscosity and a mathematical $L_2$ dissipation) considered separately. Furthermore, when the membrane and the bulk viscosity are operative simultaneously, the size of the system dictates the dominant dissipative mechanism. For vesicles smaller than a given length-scale, membrane viscosity dominates the dynamical behavior. This length-scale depends on the membrane viscosity, which can exhibit variations of several orders or magnitude depending on the amphiphiles chemistry and composition, as well as on temperature \citep{DanovDimova00}. In a significant set of situations, the regime dominated by membrane viscosity turns out to be the relevant one, as the separation length-scale is on the order of a few to tens of microns for usual lipid membranes at physiological temperatures \citep{Dimova:2006ys}, reaching hundreds of microns for liquid ordered phases \citep{Bacia:2005lr} and several millimeters for polymer vesicles \citep{Discher:1999sp,Dimova:2002ij,Zhou:2005xy}.


\begin{acknowledgments}
\noindent \emph{Acknowledgments.} We thank Laurence Grosjean and Thiwanka Wickramasooriya for their help during the early stages of this work. We gratefully acknowledge the support of the Italy-Spain Integrated Action (Grant No. HI2006-0212). MA acknowledges the support of the European Research Council under the European Community's 7th Framework Programme (FP7/2007-2013)/ERC grant agreement nr 240487, the Ministerio de Ciencia e Innovaci\'on (DPI2007-61054) and the support received through the prize ``ICREA Academia'' for excellence in research, funded by the Generalitat de Catalunya.
\end{acknowledgments}

\appendix

\section{Partial derivatives of the mean curvature}
\label{G}

Taking variations of equation (\ref{curv_e}), it follows that
\begin{equation}
\begin{split}
G[\dot{r},\dot{z}] = -\int_0^1 \kappa (H-C_0) \left(\frac{\partial H}{\partial r}\dot{r} + \frac{\partial H}{\partial r'}\dot{r}' + \frac{\partial H}{\partial r''}\dot{r}''  + \frac{\partial H}{\partial z'}\dot{z}' + \frac{\partial H}{\partial z''}\dot{z}''     \right)(2\pi a r){\rm d}u \\ 
+ \int_0^1  {\kappa} (H-C_0)^2 \left( \frac{r'}{a}\dot{r}' + \frac{z'}{a}\dot{z}' + a \dot{r} \right)\pi {\rm d}u,
\end{split}
\label{G_rz}
\end{equation}
where
\[
\frac{\partial H}{\partial r} = -\frac{z'}{a r^2},
~~~
\frac{\partial H}{\partial r'} = \frac{z''}{a^3} - \frac{3 b r'}{a^5},
~~~
\frac{\partial H}{\partial r''} = -\frac{z'}{a^3},
~~~ 
\frac{\partial H}{\partial z'} = -\frac{r''}{a^3} - \frac{3 b z'}{a^5} + \frac{1}{a r},
~~~
\frac{\partial H}{\partial z''} = \frac{r'}{a^3}.
\]

\section{Local inextensibility constraint in the Lagrangian gauge}
\label{inext}

We rewrite the axisymmetric form of the inextensibility constraint in Eq.~(\ref{inext_axi}) under the Lagrangian gauge given by Eq.~(\ref{v_rz}):
\begin{eqnarray*}
0 & = &  \frac{1}{ar} (r v_t)' - H v_n = \frac{1}{a}v_t' + \frac{r'}{ar}v_t - H v_n \\
& = & \frac{1}{a} \left[-\frac{a'}{a^2} (r'\dot{r} + z'\dot{z}) + \frac{1}{a} (r''\dot{r} + z''\dot{z} + r'\dot{r}' + z'\dot{z}') \right] + \frac{r'}{r a^2}(r'\dot{r} + z'\dot{z}) \\
&& - \frac{1}{a^2}\left( \frac{-r''z'+ r'z''}{a^2} + \frac{z'}{r}  \right) (-z'\dot{r} + r'\dot{z}) \\
& = &  \underbrace{\left[  -\frac{1}{a^4}({r'}^2 \dot{r} + r'z'\dot{z}) + \frac{1}{a^2} \dot{r}  +  \frac{1}{a^4} (-{z'}^2 \dot{r} + r'z'\dot{z}) \right]}_{0} r'' \\
&& + \underbrace{\left[    -\frac{1}{a^4}({z'}^2 \dot{z} + r'z'\dot{r}) + \frac{1}{a^2} \dot{z}  + \frac{1}{a^4} (-{r'}^2 \dot{z} + r'z'\dot{r}) \right]}_{0} z'' \\
&& + \frac{1}{a^2} (r'\dot{r}' + z'\dot{z}') +   \frac{r'}{r a^2}(r'\dot{r} + z'\dot{z}) - \frac{z'}{r a^2} (-z'\dot{r} + r'\dot{z}) \\
& = & \frac{1}{a^2} (r'\dot{r}' + z'\dot{z}') + \frac{\dot{r}}{r} = \frac{\dot{a}}{a} + \frac{\dot{r}}{r} = \frac{\partial}{\partial t}\ln ar
\end{eqnarray*}

\bibliographystyle{jfm}
\bibliography{bibdata}

\end{document}